\documentclass[11pt,a4paper]{amsart}
\usepackage{a4wide}
\usepackage{amsfonts,amsmath,amssymb,amsthm,enumerate,color,hyperref,bbm,tabularx,ifthen,twoopt,mathabx}
\usepackage{graphicx}
\usepackage[latin1]{inputenc}
\usepackage{subfig}
\usepackage{array}

\newcounter{hypA}

\def\Var{\mathop{\rm Var}\nolimits}

\def \1{\mathbbm{1}}
\def\Y{\mathbf{Y}}
\def\betabold{\boldsymbol{\beta}}
\def\X{\mathbf{X}}
\def\PX{\mathbf{P}_{\X}}
\def\A{\mathbf{A}}
\def\u{\mathbf{u}}
\def\Z{\mathbf{Z}}
\def\R{\mathbf{R}}
\def\W{\mathbf{W}}

\def\U{\mathbf{U}}

\def\e{\mathbf{e}}

\def\Id{\textrm{Id}}
\def\rset{\mathbb{R}}

\theoremstyle{remark}

\begin{document}

\title{Improving heritability estimation by a variable selection approach
in sparse high dimensional linear mixed models}

\author{A. Bonnet}
\author{C. L\'evy-Leduc}
\author{E. Gassiat}
\author{R. Toro}
\author{T. Bourgeron}
\address{Laboratoire de Mathématiques d'Orsay, Univ. Paris-Sud, CNRS, Université
Paris-Saclay, 91405 Orsay, France}
\email{elisabeth.gassiat@math.u-psud.fr}
\address{AgroParisTech/UMR INRA MIA 518}
\email{celine.levy-leduc@agroparistech.fr, anna.bonnet@agroparistech.fr}
\address{Human Genetics and Cognitive Functions,Institut Pasteur, Paris, France}
\email{rto@pasteur.fr, thomasb@pasteur.fr}

\maketitle




\begin{abstract}
Motivated by applications in neuroanatomy, we propose a novel methodology for estimating the heritability which corresponds to the
proportion of phenotypic variance which can be explained by genetic factors. Estimating this quantity for 
neuroanatomical features is a fundamental challenge in psychiatric disease research. Since the phenotypic variations may only be due to
a small fraction of the available genetic information, we propose an estimator of 
the heritability that can be used in high dimensional sparse linear 
mixed models. Our method consists of three steps. Firstly, a variable selection stage is performed in order to recover the support 
of the genetic effects  -- also called causal variants -- that is to find the genetic effects which really explain the phenotypic variations.
Secondly, we propose a maximum likelihood strategy for estimating the heritability which only takes into account
 the causal genetic effects found in the first 
step. Thirdly, we compute the standard error and the 95\% confidence interval associated to our heritability estimator thanks to a nonparametric
bootstrap approach. Our contribution consists in providing an estimation of the heritability with standard errors 
substantially smaller than methods without variable selection
\textcolor{black}{when the genetic effects are very sparse. Since the
real genetic architecture is in general unknown in practice, we also propose an empirical criterion which allows the user to decide
whether it is relevant to apply a variable selection based approach or not.} We illustrate the performance
of our methodology on synthetic and real neuroanatomic data coming from the Imagen project. We also show that our approach 
has a very low computational burden and is very efficient from a statistical point of view.
\end{abstract}

\section{Introduction}









For many complex traits in human population, there exists a huge gap between the genetic variance
explained by population studies and the variance explained by specific variants found thanks to
genome wide association studies (GWAS). This gap has been called by \cite{maher:2008} and \cite{manolio:2009}
the ``dark matter'' of the genome or the ``dark matter'' of heritability. 
Various population studies have shown that up to $80\%$ of the variability of neuroanatomical phenotypes such as the brain volume could be explained by genetic factors, see for instance \cite{stein:2012}. This result is very important since several psychiatric disorders are shown to be associated to neuroanatomical changes, for instance macrocephaly and autism \cite{steen:2006} or reduced hippocampus and schizophrenia \cite{amaral:2008}. Estimating properly the impact of the genetic background on neuroanatomical changes is a crucial challenge in order to determine afterwards if this background can either be a risk factor or a protective factor from developing psychiatric disorders. The GWAS studies performed for instance by \cite{stein:2012} identified genetic variants involved in the neuroanatomical diversity, which contributes to understand the impact of genetic factors. However, in the course of these studies, it is shown that this approach only explains a small proportion of the phenotypic variance. In order to understand the nature of the genetic factors responsible for major variations of the brain volume, \cite{toro:2014} used linear mixed models (LMM) to consider the effects of all the common genetic diversity characterized by the Single Nucleotide Polymorphisms (SNPs). This approach had been suggested by \cite{yang:lee:goddard:visscher:2011} to study the effects of the SNPs on the height variations.
The model they considered is a LMM defined as follows:
\begin{equation}\label{eq:modele}
\Y=\X\betabold+\Z\u+\e\;, 
\end{equation}
where $\Y=(Y_1,\dots,Y_n)'$ is the vector of observations (phenotypes), $\X$ is a $n\times p$ matrix of predictors, $\betabold$
is a $p\times 1$ vector containing the unknown linear effects of the predictors, $\Z$ is the genetic information matrix, $\u$ and $\e$ correspond to the random
effects. More precisely, $\Z$ is a version of $\W$ with centered and normalized columns, where $\W$ is defined as follows: $W_{i,j}=0$ (resp. 1, resp. 2) if the genotype of the $i$th individual at locus $j$ is $qq$
(resp. Qq, resp. QQ) where $p_j$ denotes the frequency of the allele q at locus $j$.
In (\ref{eq:modele}), the vector $\e$ corresponds to the environment effects and the vector $\u$ corresponds to the genetic random effect, that is the $j$-th component of $u$ is the effect of the $j$-th SNP on the phenotype. 
In the modeling of \cite{yang:lee:goddard:visscher:2011}, all the SNPs have an effect on the considered phenotype, that is
\begin{equation}\label{eq:distrib_u_e_yang}
\u\sim\mathcal{N}\left(0,{\sigma_u^\star}^2\Id_{\rset^n}\right) \textrm{ and } \e\sim\mathcal{N}\left(0,{\sigma_e^\star}^2\Id_{\rset^n}\right).
\end{equation}
The \textcolor{black}{covariance matrix} of $\Y$ can thus be written as:
$$
\Var(\Y)=N{\sigma_u^\star}^2 \R +{\sigma_e^\star}^2 \Id_{\rset^n}\;,\textrm{ where }\R=\frac{\Z\Z'}{N}\;,
$$
and the parameter $\eta^\star$ defined as
\begin{equation}
\label{eq:herit_yang}
 \eta^\star
=\frac{N{\sigma_u^\star}^2}{N{\sigma_u^\star}^2+{\sigma_e^\star}^2}\;
\end{equation}
is commonly called the heritability (\cite{yang:lee:goddard:visscher:2011},\cite{pirinen:donnelly:spencer:2013}), and corresponds to the proportion of phenotypic variance which is determined by all the SNPs.
 
 Since all SNPs are not necessarily causal, it seems more realistic to extend the previous modeling by assuming that the genetic random effects can be sparse, that is only a proportion $q$ of the components of $\u$ are non null:
\begin{equation}\label{eq:distrib_u_e_comp_nulles}
u_i\stackrel{i.i.d.}{\sim}(1-q)\delta_0+q\mathcal{N}(0,{\sigma_u^\star}^2),\textrm{ for all } 1\leq i\leq N,
\end{equation}
where $q$ is in $(0,1]$, and $\delta_0$ is the point mass at $0$.
Then the definition of $\eta^\star$ has to be adjusted as follows:
\begin{equation}
 \eta^\star
=\frac{Nq{\sigma_u^\star}^2}{Nq{\sigma_u^\star}^2+{\sigma_e^\star}^2}\;.
\end{equation}
It corresponds to the proportion of phenotypic variance which is due to a certain number of causal SNPs which are, obviously, unknown.
Let us emphasize that, in most applications, the proportion $q$ of causal SNPs is also unknown, and that it may happen that the scientist has no idea how small $q$ is.

When $q=1$, that is when considering the modeling \eqref{eq:distrib_u_e_yang}, most proposed approaches to estimate the heritability  derive from a likelihood methodology. We can quote for instance the REstricted Maximum Likelihood (REML) strategies, originally proposed by \cite{patterson:thompson:1971} and then developed in \cite{searle:casella:mcculloch:1992}. Several approximations of the REML algorithm have also been proposed, see for instance the software EMMA proposed by \cite{pirinen:donnelly:spencer:2013} or the software GCTA (\cite{yang:lee:goddard:visscher:2011},\cite{Yang:2010}). 

 We proposed in \cite{nous} another method based on a maximum likelihood strategy to estimate the heritability and implemented in the R package HiLMM. We proved in \cite{nous} the 
following theoretical result: though the computation of the likelihood is based on the modeling assumption \eqref{eq:distrib_u_e_yang}, the estimator is consistent (unbiased) under the less restrictive modeling assumption (\ref{eq:distrib_u_e_comp_nulles}). We believe this consistency result  
remains true for the estimators produced using the algorithms REML, EMMA, GCTA. But we also proved that, when $q\neq 1$, the standard error is not the one computed by the softwares when $q=1$ and  may be very large. We obtained a theoretical formula for the asymptotic variance of the estimator (depending in particular on $q$) and conducted several numerical experiments to understand how this asymptotic variance gets larger depending on the various quantities,  in particular with respect to $q$ and the ratio $n/N$. We observed that this variance indeed gets larger when $q$ gets smaller, so that the accuracy of the heritability estimator is slightly deteriorated when all SNPs are not causal. Thus,  a first problem  is to find a method able to produce an estimator with smaller standard error than those obtained using only likelihood strategies.
 Also, since this standard error depends on $q$, a second problem 
 is to produce a confidence interval one could trust without knowing $q$.
 
 The goal of this paper is to address both problems. 
 The results we obtained in  \cite{nous} suggest the following. If we knew the set of causal SNPs, then, considering only this (small) subset in the genetic information matrix, we would obtain with HiLMM an estimator having a smaller standard error than when using all SNPs in the genetic information matrix. Thus, our new practical method  contains a variable selection step.
 
Variable selection and signal detection  in high dimensional linear models have been extensively studied in the past decade and there are many papers on this subject. 
 Among them, we can  quote 
\cite{meinshausen:buhlmann:2010} and \cite{blanchard:2014} about variable selection and  references therein. 
The case of high dimensional 
 mixed models has received little attention.
As far as variable selection methods in the random effects of LMM are concerned, we are only aware of the work of \cite{Fan_Li} and \cite{bondell_2010}. Let us mention that regarding the estimation of heritability with possible sparse effects, there is also the bayesian approach of \cite{guan:stephens:2011} and \cite{zhou:carbonetto:stephens:2013}, which proposes an interesting estimator for the heritability but which is computationally very demanding.
Notice that, in our framework, we are not far from the situation
 for which  it is proved in \cite{verzelen:2012} that the support cannot be fully recovered, 
which happens 
 when  $Nq\log(1/q) >> n$.
The variable selection step we propose takes elements from both ultrahigh dimension methods (\cite{fan:lv:2008}, \cite{Ji_Jin}, \cite{meinshausen:buhlmann:2010}) and classical variable selection techniques (\cite{Tibshirani:96}). 


The second step of our method is to apply HiLMM using the selected
subset of causal SNPs produced by the first step. Finally, we propose
a non parametric bootstrap procedure to get confidence intervals with
prescribed coverage. The whole procedure requires only a few minutes
of computation.

To conclude,
we propose in this paper a very fast method to estimate the
heritability and construct a confidence interval substantially smaller
than without variable selection \textcolor{black}{when the genetic effects are very sparse.
Since the real genetic architecture is in general unknown in practice, we also propose an empirical criterion which allows the user to decide
whether it is relevant to apply a variable selection based approach or not.}
Our method has also the advantage to return a list of SNPs possibly
involved in the variations of a given quantitative feature. 
This set of SNPs can further be analyzed from a biological point of view.

The paper is organized as follows. Section \ref{sec:data_set}
describes the data set which motivated our work. Section
\ref{sec:method} provides the detailed description of the method, and
Section \ref{sec:num_study} displays the results of the numerical
study. They were obtained by using the R package EstHer that we
developed and which is available from the Comprehensive R Archive
Network (CRAN). The simulation results illustrate the performance of
our method on simulations and show that it is very efficient from a
statistical point of view. \textcolor{black}{In Section
  \ref{sec:criterion}, we provide an empirical criterion to help the
  user to decide whether it is relevant to apply a variable selection
  based approach or not. In Section \ref{sec:comparison}, we propose a
thorough comparison of our approach with other methods in terms of statistical and
numerical performances.}
Finally, the results obtained on the brain data described in Section \ref{sec:data_set} can be found in Section \ref{sec:real_data}.
We also provide a discussion section at the end of the paper.


 \section{Description of the data}
 \label{sec:data_set}

\begin{figure}[h!]
 \begin{center}
 \includegraphics[scale=0.4,trim=0mm 110mm 160mm 20mm,clip]{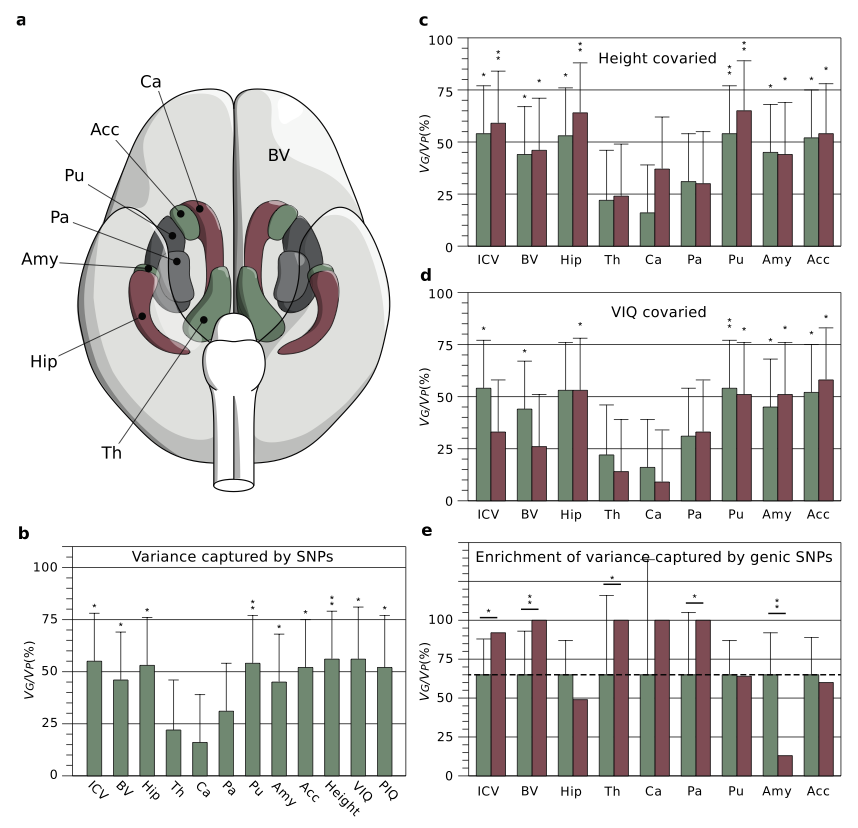}
 \caption{Different regions of the brain (this figure is taken from \cite{toro:2014}).}
 \label{fig:cerveau}
    \end{center}
 \end{figure}

We worked on data sets provided by the European project Imagen, which is a major study on mental health and risk taking behaviour in teenagers. The research program includes questionnaires, interviews, behaviour tests, neuroimaging of the brain and genetic analyses.
We will focus here on the genetic information collected on
approximately $2000$ teenagers as well as measurements of the volume
of several features: the intracranial brain volume (icv), the thalamus
(th), the caudate nucleus (ca), the amygdala (amy), the globus
pallidus (pa), the putamen (pu), the hippocampus (hip), the nucleus
accubens (acc) and the total brain volume (bv). Figure
\ref{fig:cerveau}, which comes from \cite{toro:2014}, is a schematic
representation of these different areas of the brain. 
The data set contains $n=2087$ individuals and $N=273926$ SNPs, as well as a set of fixed effects, which in our case are the age (between 12 and 17), the gender and the city of residency (London, Nottingham, Dublin, Dresden, Berlin, Hamburg, Mannheim and Paris).

 In the following, our goal will thus be to provide a method to estimate the heritability of these neuroanatomical features.

\section{Description of the method}\label{sec:method}

The method that we propose can be split into two main parts: the first one consists in a variable selection approach
and the second one provides an estimation of the heritability and the associated $95\%$ confidence interval which is computed by using 
non parametric bootstrap.

At the beginning of this section we shall consider the case where there is no fixed effects, that is 
\begin{equation}\label{eq:modele_without_fixed_effects}
\Y=\Z\u+\e
\end{equation}
but we explain at the end of this section how to deal with fixed effects. Let us first describe our variable selection 
method which consists of three steps.


\subsection{Variable selection}\label{sec:var_select}

Inspired by the ideas of \cite{fan:lv:2008}, we do not directly apply a Lasso type approach since we are
in an ultra-high dimension framework. Hence, we start our variable selection stage
by the SIS (Sure Independence Screening) approach, as suggested by \cite{fan:lv:2008}, 
in order to select the components of $\u$ which are the most correlated to the
response $\Y$ and then we apply a Lasso criterion which depends on a regularization parameter $\lambda$. This regularization parameter
is usually chosen by cross validation but here we decided to use the stability selection approach devised by \cite{meinshausen:buhlmann:2010}
which provided better results in our framework.

\subsubsection*{Step 1: Empirical correlation computation}

The first step consists in reducing the number of relevant columns of $\Z$ by trying to remove those associated to null
components in the vector $\u$. For this, we use the SIS (Sure Independence Screening) approach proposed by \cite{fan:lv:2008}
and improved by \cite{Ji_Jin} in the ultra-high dimensional framework. More precisely, we compute for each column $j$ of
$\Z$:
$$C_j=\left\vert \sum Y_i Z_{i,j}\right\vert , $$
and we only keep the $N_{\textrm{max}}$ columns of $\Z$ having the largest $C_j$. In practice,  we choose the conservative value $N_{\textrm{max}}=n$, inspired by the comments of \cite{fan:lv:2008} on the choice of $N_{\textrm{max}}$. 

In the sequel, we denote by $\Z_{\textrm{red}}$ the matrix containing these $n$ relevant columns.
This first step is essential for our method. Indeed, on the one hand, it substantially decreases the computational burden of our approach
and on the other hand, it reduces the size of the data and thus makes classical variable selection tools efficient.

\subsubsection*{Step 2: LASSO criterion \textcolor{black}{and stability selection}}

In order to refine the set of columns (or components of $\u$) selected in the first step and to remove the
remaining null components in the vector $\u$, we apply a Lasso criterion originally devised by
\cite{Tibshirani:96} which has been used in many different contexts and has been thouroughly theoretically studied. 
It consists in minimizing with respect to $u$ the following criterion:
\begin{equation}\label{eq:lasso}
\textrm{Crit}_\lambda(u)=\| \Y - \Z_{\textrm{red}} u \|_2^{2} + \lambda \| u \|_1\;,
\end{equation}
which depends on the parameter $\lambda$ and where $\|x\|_2^2=\sum_{i=1}^p x_i^2$ and $\|x\|_1=\sum_{i=1}^p |x_i|$  for
$x=(x_1,\dots,x_p)$. The choice of the regularization parameter $\lambda$ is crucial since its value may strongly affect
the selected variables set. Different approaches have been proposed for choosing this parameter such as cross-validation 
which is implemented for instance in the \verb|glmnet| R package. 
\textcolor{black}{Here we shall use the following strategy based on the stability selection
 proposed by \cite{meinshausen:buhlmann:2010}.}

%
%
%
%
%
The vector of observations $\Y$ is randomly split into several subsamples of size $n/2$. 
\textcolor{black}{For each subsample, we apply the LASSO criterion for a fixed parameter $\lambda$
and the selected variables are stored}. 
Then, for a given threshold, we keep in the final
set of selected variables only the variables appearing a number of times larger than this threshold. 
In practice, we generated $50$ subsamples of $\Y$ \textcolor{black}{and we chose the parameter $\lambda$ as the
smallest value of the regularization path. As explained in \cite{stab_selec}, such a choice of $\lambda$ ensures that some overfitting
occurs and hence that the set of selected variables is large enough to include the true variables with high probability.}

The matrix $\Z$ containing only the final set of selected columns will be denoted by $\Z_{\textrm{final}}$ in the following, where 
$N_{\textrm{final}}$ denotes its number of columns.

The threshold has to be chosen  carefully: keeping too many columns in  $\Z_{\textrm{final}}$ could indeed lead
to overestimating the heritability and, on the contrary, removing too many columns of $\Z$ could
lead to underestimating the heritability. In the ``small $q$'' situations where
it is relevant to use a variable selection approach a range of thresholds in which the heritability estimation is stable will appear
as suggested by \cite{stab_selec}.
In practice, we simulate observations $\Y$ satisfying (\ref{eq:modele_without_fixed_effects}),
by using the matrix $\Z$, for different values of $q$ and for different values $\eta^\star$ and we observe that this stability region
for the threshold appear for small values of $q$. This procedure is further illustrated in Section \ref{sec:num_study}.


\subsection{Heritability estimation and confidence interval}

\subsubsection{Heritability estimation}
\label{subsec:herit_estim}
For estimating the heritability, we used the approach that we proposed in \cite{nous}. 
It is based on a maximum likelihood strategy and was implemented in the R package HiLMM. Let us recall how this method works.

In the case where $q=1$, which corresponds to the non sparse case,
$$
\Y\sim\mathcal{N}\left(0,\eta^\star{\sigma^\star}^2 \R+(1-\eta^\star){\sigma^\star}^2 \Id_{\rset^n}\right),$$ with  $\sigma^{\star 2}=N \sigma_u^{\star 2} + \sigma_e^{\star 2}$ and 
$\R=\Z_{\textrm{final}}\Z_{\textrm{final}}'/N_{\textrm{final}},$ 
\textcolor{black}{where $\Z_{\textrm{final}}$ denotes the matrix $\Z$ in which the columns selected in the variable selection step 
described in Section \ref{sec:var_select} are kept.}

Let $\U$ be defined as follows: $\U'\U=\U\U'=\Id_{\rset^n}$ and
$\U\R\U'=\textrm{diag}(\lambda_1,\dots,\lambda_n)$, where the last quantity denotes the diagonal matrix having its diagonal entries
equal to $\lambda_1,\dots,\lambda_n$.
Hence, in the case where $q=1$,
\textcolor{black}{
\begin{multline}\label{eq:widetildeY}
\widetilde{\Y}=\U'\Y \sim\mathcal{N}(0,\Gamma)\\
\textrm{ with } \Gamma=\textrm{diag}(\eta^\star{\sigma^\star}^2\lambda_1+(1-\eta^\star){\sigma^\star}^2
,\dots,\eta^\star{\sigma^\star}^2\lambda_n+(1-\eta^\star){\sigma^\star}^2),
\end{multline}
}
where the $\lambda_i$'s are the eigenvalues of $\R$.

We propose to define $\hat{\eta}$ as a
maximizer of the log-likelihood
\begin{equation}\label{eq:Ln}
L_n(\eta)= -\log \left(\frac{1}{n} \sum_{i=1}^{n} \frac{\widetilde{Y}_i^{2}}{\eta(\lambda_i-1)+1}\right) 
-\frac{1}{n}\sum_{i=1}^{n} \log\left(\eta(\lambda_i-1)+1\right)\;,
\end{equation}
where the $\widetilde{Y}_i$'s are the components of the vector $\widetilde{\Y}=\U'\Y$. 


We now explain how to obtain accurate confidence intervals for the heritability by using a non parametric bootstrap approach.

\subsubsection{Bootstrap confidence interval}

We used the following procedure:
\textcolor{black}{
\begin{itemize}
\item[-] Step 1: We estimate $\eta^{\star}$ and ${\sigma^{\star}}^2$ by using our approach described in the previous subsection.
The corresponding estimators are denoted $\hat{\eta}$ and $\hat{\sigma}$.
\item[-] Step 2: 
We compute $\Y_{\textrm{new}}=\hat{\Gamma}^{-1/2}\widetilde{\Y}$, where \textcolor{black}{$\widetilde{\Y}$ is defined
in (\ref{eq:widetildeY}) and $\hat{\Gamma}$ has the same structure as
$\Gamma$ defined in (\ref{eq:widetildeY})} except that $\eta^\star$ and $\sigma^\star$ are replaced 
by their estimators $\hat{\eta}$ and $\hat{\sigma}$, respectively. 
\item[-] Step 3: We create $K$ vectors $(\Y_{\textrm{new},i})_{1\leq i\leq K}$ 
from $\Y_{\textrm{new}}$ by randomly choosing each of its components among those of $\Y_{\textrm{new}}$.
\item[-] Step 4: We then build $K$ new vectors $(\widetilde{\Y}_{\textrm{samp},i})_{1\leq i\leq K}$ as follows: 
$\widetilde{\Y}_{\textrm{samp},i}=\hat{\Gamma}\Y_{\textrm{new},i}$. For each of them we estimate the heritability.
We thus obtain a vector \textcolor{black}{of heritability estimators} $(\hat{\eta}_1,...,\hat{\eta}_K)$.
\item[-] Step 5: For obtaining a 95\% bootstrap confidence interval, we order these values of $\hat{\eta}_k$ and keep
the ones corresponding to the $\lfloor 0.975\times K\rfloor$ largest and the $\lfloor 0.025\times K \rfloor$ smallest, where $\lfloor x \rfloor$ denotes the integer part of $x$. These values define
the upper and lower bounds of the 95\% bootstrap confidence interval for the heritability $\eta^\star$, respectively.
\end{itemize}
}

A bootstrap estimator of the variance can be obtained by computing the empirical variance estimator of the 
$\hat{\eta}_k$'s. In practice, we chose $K=80$ replications.


\textcolor{black}{In Step 2 of the previous algorithm, we should be in the non sparse case $q=1$ thanks to the variable selection stage.
Hence, the covariance matrix of $\Y_{\textrm{new}}$ should be close to identity.
}

\textcolor{black}{Observe that our resampling technique is close to the one proposed by \cite{GEPI:GEPI21893} for building permutation tests
in linear mixed models.
}
\subsection{Additional fixed effects}

The method described above does not take into account the presence of fixed effects. For dealing with such effects we propose to use the 
following method, which mainly consists in projecting the observations onto the orthogonal of $\mathrm{Im}(\X)$, the image of $\X$, 
to get rid of the fixed effects. 
In practice, instead of considering $\Y$ and $\Z$ we consider $\tilde{\Y}=\A'\Y$ and $\tilde{\Z}=\A'\Z$, where $A$ is a $n \times (n-d)$ matrix 
($d$ being the rank of the fixed effects matrix), such that $\A\A'=\PX$, $\A'\A=\textrm{Id}_{\rset^{n-d}}$ and $\PX= \textrm{Id}_{\rset^n} - \X(\X'\X)^{-1}\X'$. 
This procedure was for instance used by \cite{Fan_Li}.

\section{Numerical study}
 \label{sec:num_study}

We present in this section the numerical results obtained with our approach which is implemented in the R package EstHer.

\subsection{Simulation process}
\label{subsec:simul}
Since in genetic applications, the number $n$ of individuals is very small with respect to the number $N$ of SNPs, we chose $n=2000$ and $N=100000$ 
in our numerical study. We also set $\sigma_u^{\star2}=1$, \textcolor{black}{we shall
consider different values for $q$ and we shall
change the value of $\sigma_e^{\star}$ in order to have the following
values for $\eta^\star$: 0.4, 0.5, 0.6 and 0.7.} 
We generate a matrix $\W$ such that its columns $W_j$ are independent binomial random variables of parameters $n$ and $p_j$, where $p_j$ 
 is randomly chosen in $[0.1,0.5]$. We compute $\Z$ 
by centering and empirically normalizing the matrix $\W$.
The random effects are generated according to Equation (\ref{eq:distrib_u_e_comp_nulles}) and then we compute a vector of observations such 
that $\Y=\Z \u+\e$.

We can make two important comments about the previous simulation process. Firstly, we generated a matrix $\W$ with independent columns, 
that is we assume that the SNPs are not correlated. Since this assumption may not be very realistic in practice, we provide in Section 
\ref{subsec:vraie_Z} some additional simulations where the generated matrix $\W$ has been replaced by the real matrix $\W$ coming from 
the IMAGEN project. Secondly, we did not include fixed effects but we show some results in Section \ref{subsec:fixed_effects}
when fixed effects are taken into account. 


\subsection{Results \textcolor{black}{in very sparse scenarios}}
\label{subsec:res_q_small}

\textcolor{black}{In this section, we shall focus on the performances of
  our method in a very sparse scenario, that is 100 causal SNPs out of
  100,000. We will describe all the results in terms of heritability
  estimation, support recovery and computational times 
in this particular case, then we will study other sparsity scenarios.}

\subsubsection{Choice of the threshold}
\label{subsec:thresh_choice}
\begin{figure}[!htbp]
  \centering
   \includegraphics[trim=10mm 15mm 0mm 20mm,clip,width=.45\textwidth]{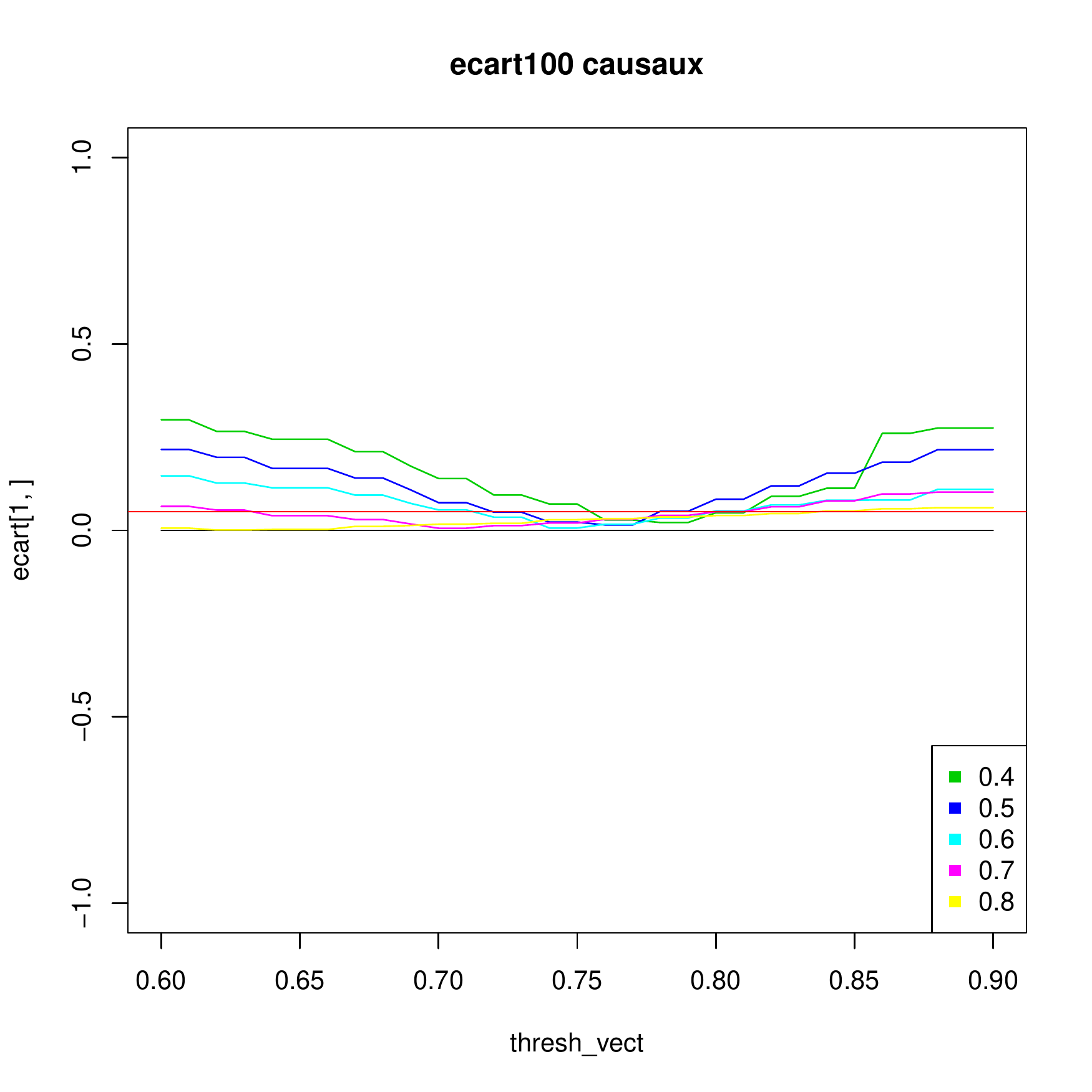} 
  \caption{Absolute difference between $\eta^\star$ and $\hat{\eta}$ for thresholds from 0.6 to 0.9 and for 
$q= 10^{-3}$ (100 causal SNPs).
}
  \label{fig:choix_seuil}
\end{figure}

In order to determine the threshold, we apply the procedure described in Section \ref{sec:var_select} 
and \ref{subsec:herit_estim}.
Figure \ref{fig:choix_seuil} 
displays the mean of the absolute value of the difference
between $\eta^\star$ and the estimated value $\hat{\eta}$ for different thresholds and for different values of $\eta^\star$ obtained from
10 replications. 
We can see from this figure that in the case where the number of causal SNPs is relatively small: 100, that is 
$q=10^{-3}$, our estimation procedure provides relevant estimations of the heritability for a range of thresholds around 0.75.
Moreover, the optimal threshold leading to the smallest gap between $\hat{\eta}$ for different values of $\eta^\star$
is 0.76.
\textcolor{black}{We will use this value in the following numerical study. However, the way of choosing the threshold will be further discussed, especially in the section dedicated to the study of the genetic data.}

\subsubsection{Confidence intervals}
 We use the non parametric boostrap approach described in Section \ref{sec:method} in order to compute the confidence intervals associated to the estimations of the heritability. Table \ref{tab:table1} shows that the $95\%$ confidence intervals obtained by bootstrap and the empirical confidence intervals are very similar.
\textcolor{black}{The empirical confidence intervals are computed as follows: the different estimations of $\eta^\star$ obtained along
the different replications are ordered, the
$\lfloor 0.975\times M\rfloor$ largest and the $\lfloor 0.025\times M\rfloor$ smallest values correspond to the upper (resp. lower) bound of the 95\% empirical confidence interval. Here, $\lfloor x\rfloor$ denotes the integer part of $x$ and $M$ is the number of replications.
From Table \ref{tab:table1}, we can see that the empirical confidence intervals are included in the bootstrap intervals, which means that
our approach provides conservative intervals.}

\begin{table}
   \caption{ 95 \% confidence intervals for $\hat{\eta}$ obtained empirically and by our Bootstrap method.}
\begin{tabular}{|l| l |l |l| l|}
  \hline
 $ \eta^{\star} $ & 0.4 & 0.5 & 0.6 & 0.7  \\
  \hline
 Bootstrap & [0.353 ; 0.503] & [0.413 ; 0.565] & [0.494 ; 0.654] & [0.596 ; 0.738] \\
  \hline   
  Empirical & [0.391 ; 0.470] & [0.449 ; 0.542] & [0.496 ; 0.645] & [0.618 ; 0.720]   \\
  \hline
\end{tabular}
\label{tab:table1}
\end{table}

\subsubsection{Comparison \textcolor{black}{between the methods with and without selection}}

\begin{figure}[!ht]
\begin{tabular}{cc}
  \includegraphics[trim=0mm 15mm 0mm 20mm,clip,width=.4\textwidth] {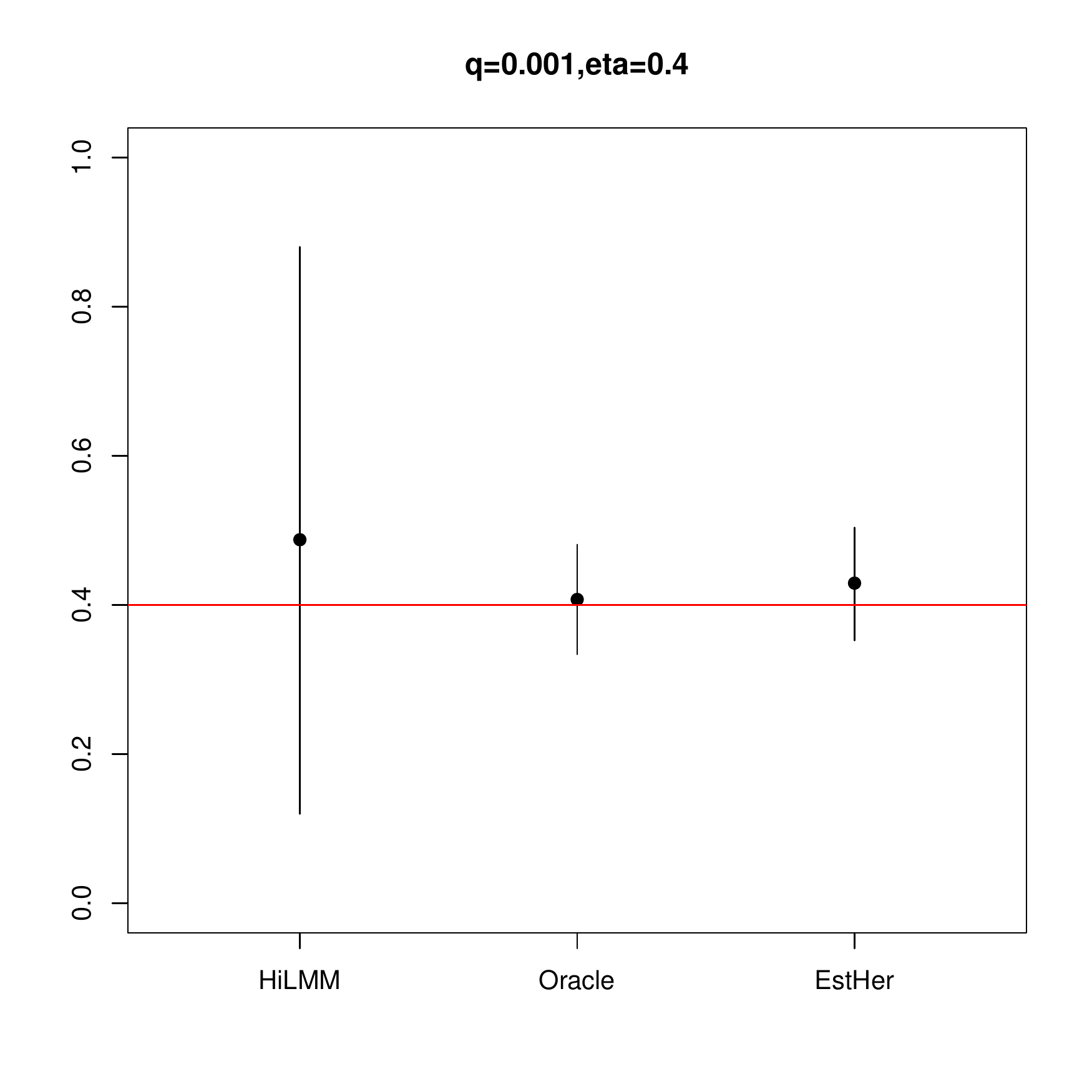}
  & \includegraphics[trim=0mm 15mm 0mm 20mm,clip,width=.4\textwidth]{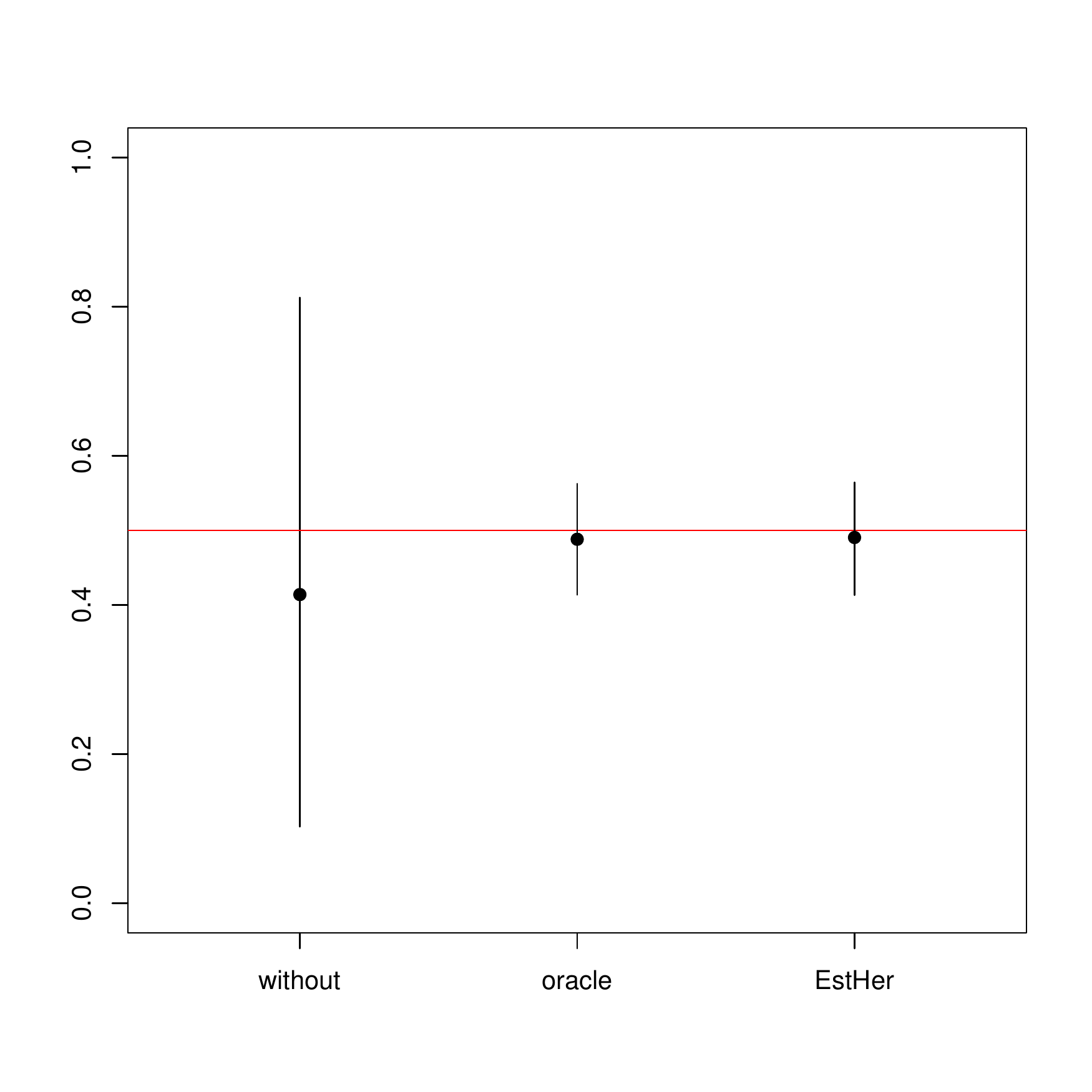} \\
(a) & (b) \\
\includegraphics[trim=0mm 15mm 0mm 20mm,clip,width=.4\textwidth]{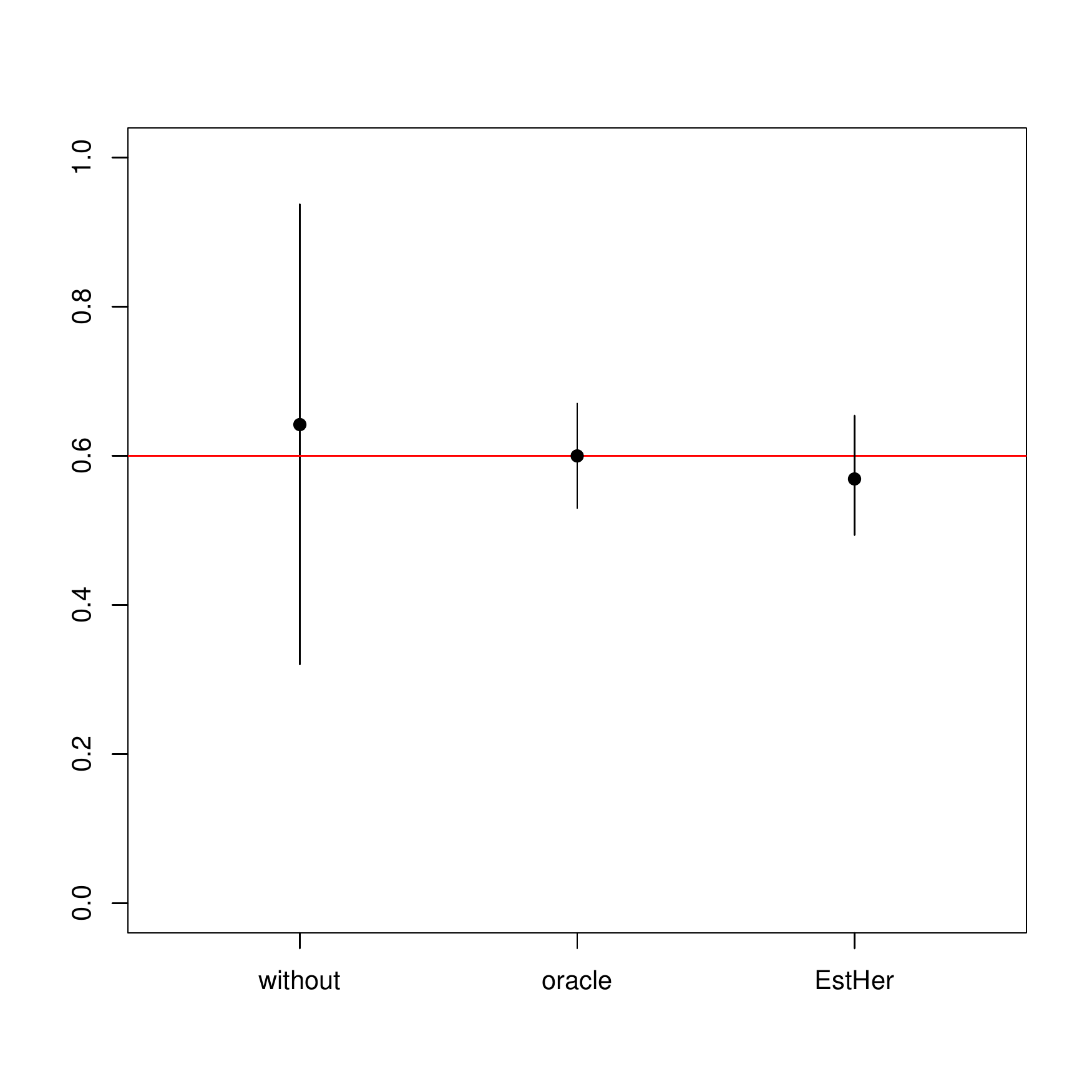} 
& \includegraphics[trim=0mm 15mm 0mm 20mm,clip,width=.4\textwidth]{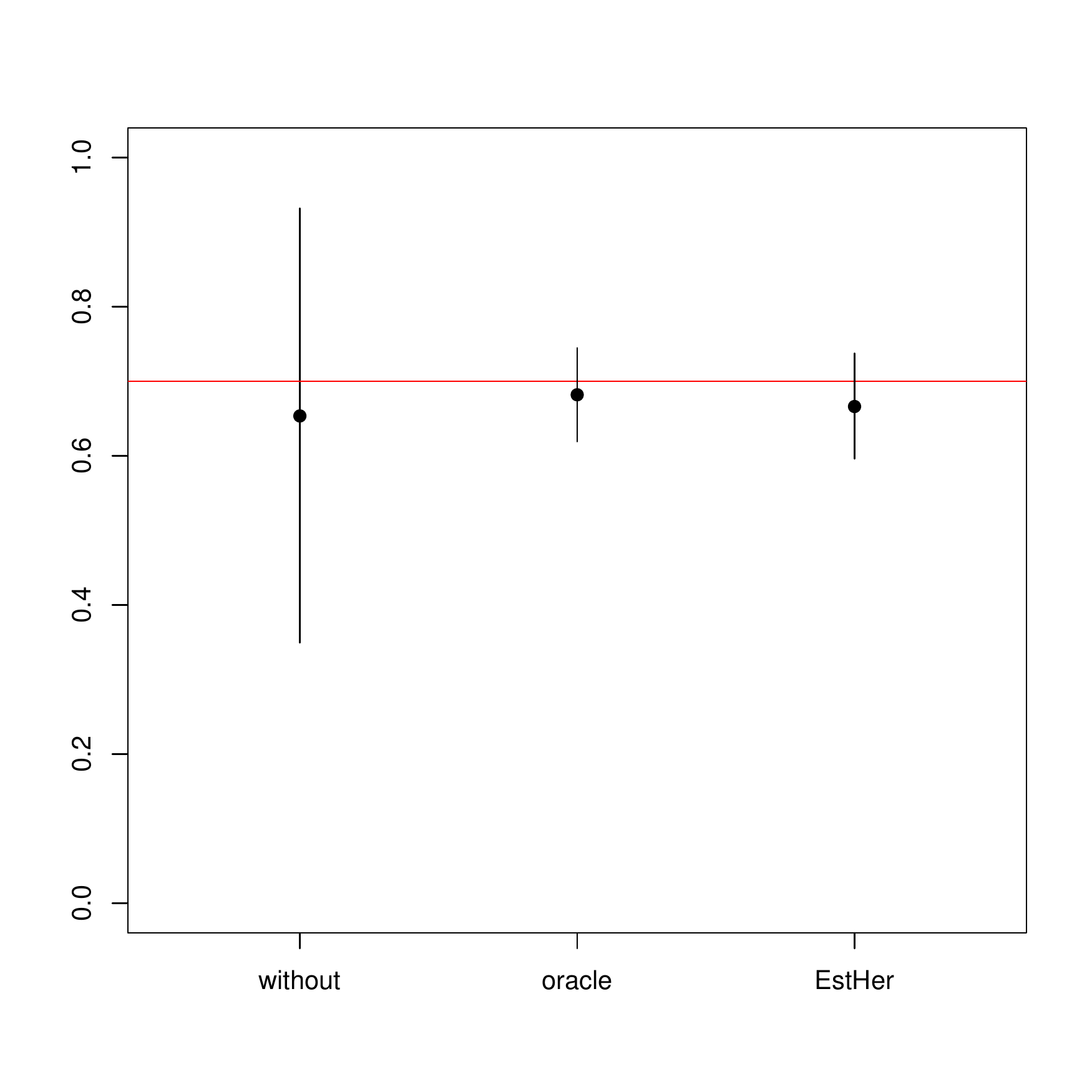} \\
(c) & (d)\\
\end{tabular}
\caption{Estimation of the heritability and the corresponding 95\% confidence intervals when $q$ =$10^{-3}$, and for different values of $\eta^{\star}$ : \textcolor{black}{(a) $\eta^{\star}=0.4$,} (b) $\eta^{\star}=0.5$, (c) $\eta^{\star}=0.6$,
(d) $\eta^{\star}=0.7$. The means of the heritability estimators (displayed with black dots), the means of the lower and upper bounds of the 95\%
confidence intervals are obtained from 20 replicated data sets for the
different methods: without selection (``without''), ``oracle'' which knows the position of the null components
 and EstHer. The horizontal gray line corresponds to the value of $\eta^{\star}$.}
  \label{fig:seuil_076}
\end{figure}

Our results are compared to those obtained if we do not perform the selection before the estimation, \textcolor{black}{that is with the method
implemented in HiLMM (''without'')}, but also with an approach which assumes the position of the non null components to be known (oracle). 
The results are displayed in Figure \ref{fig:seuil_076} and in Table \ref{tab:table2}. In this table, the confidence intervals displayed
for the lines ''Oracle'' and ''without'' are obtained by using the asymptotic variance derived in \cite{nous} which corresponds
to the classical inverse of the Fisher information in the case $q=1$. We observe that our method without the selection step provides similar results, that is almost no bias but a very large variance due 
to the framework $N\gg n$. 
Our method EstHer considerably reduces the variance compared to this method and exhibits performances close to those of the oracle 
approach which, contrary to our approach, knows the position of the non null components.

\begin{table}
   \caption{95 \% confidence intervals for $\hat{\eta}$ obtained by our approach, GCTA, the oracle approach and the approach without selection 
(``without'').}
\begin{tabular}{|l| l |l |l| l|}
  \hline
 $ \eta^{\star} $ & 0.4 & 0.5 & 0.6 & 0.7  \\
  \hline
  EstHer & [0.353 ; 0.503] & [0.413 ; 0.565] & [0.494 ; 0.654] & [0.596 ; 0.738] \\
  \hline   
  Oracle & [0.362 ; 0.472] & [0.414 ; 0.563] & [0.529 ; 0.670] & [0.619 ; 0.745] \\
  \hline
  without & [0.120 ; 0.880] & [0.102 ; 0.812] & [0.320 ; 0.938] & [0.349 ; 0.932] \\
  \hline
\end{tabular}
\label{tab:table2}
\end{table}

\subsubsection{Additional fixed effects}
\label{subsec:fixed_effects}

We generated some synthetic data according to the process described in Section \ref{subsec:simul} but we added a matrix of fixed effects 
containing two colums. Figure \ref{fig:fix_vrai} (a) displays the corresponding results which show that the presence of fixed effects does not alter the heritability estimation.

\subsubsection{Simulations with the matrix $\W$ of the IMAGEN data set}
\label{subsec:vraie_Z}
\begin{figure}[!ht]
  \centering
\begin{tabular}{cc}
  \includegraphics[width=.45\textwidth] {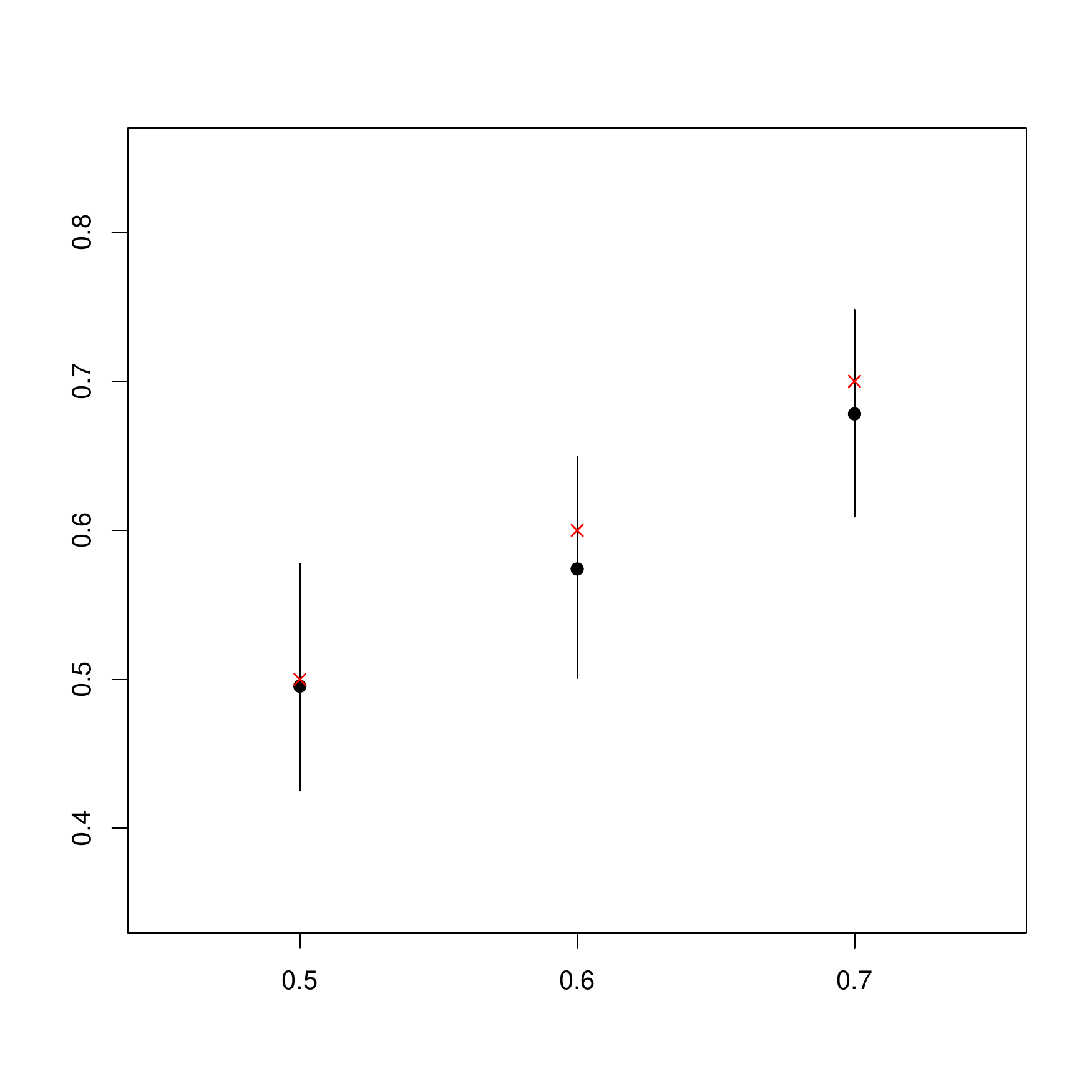}&
  \includegraphics[width=.45\textwidth]{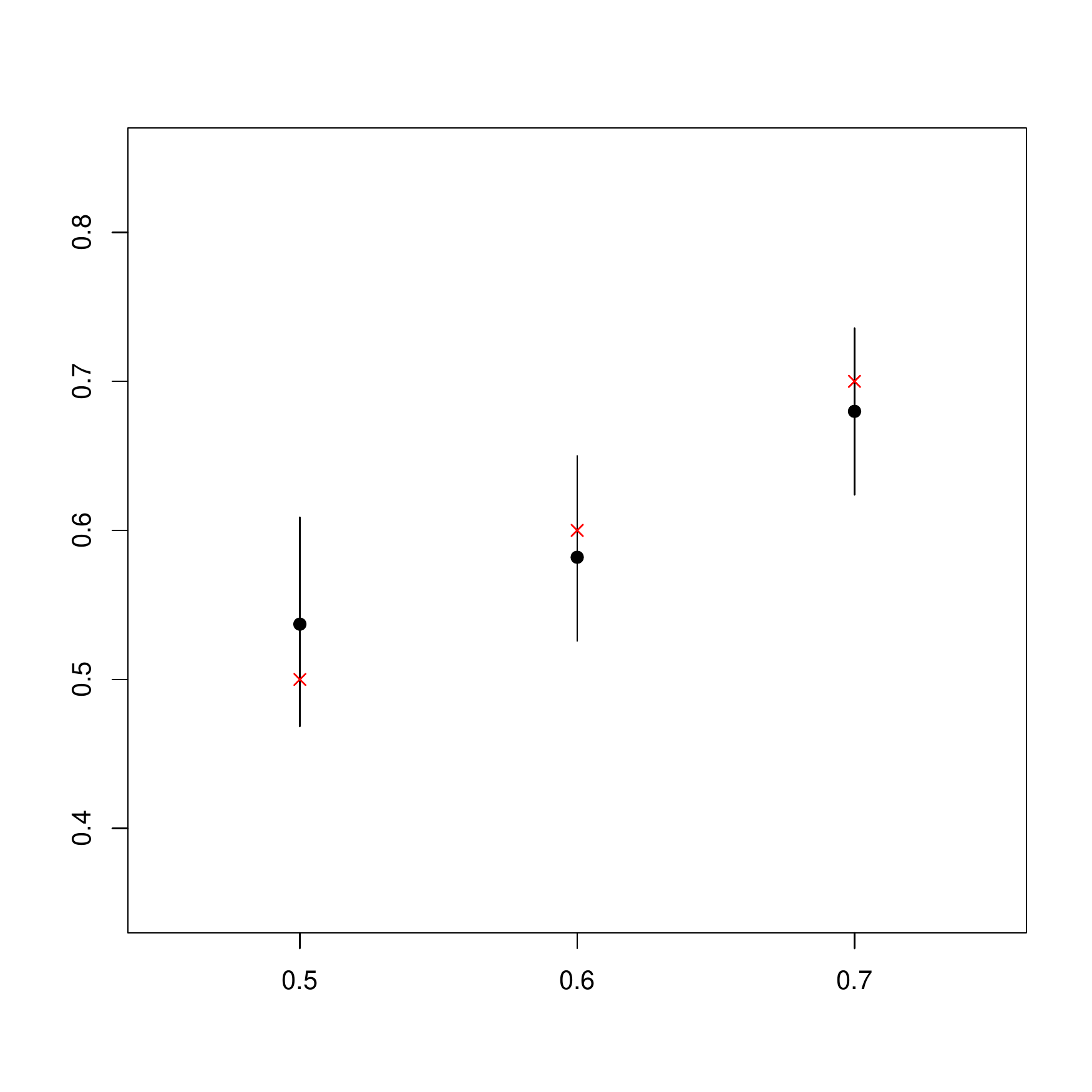}\\
(a) & (b)\\
\end{tabular}
\caption{Estimated value of the heritability with 95 \% confidence intervals. The results are displayed for several values of $\eta^\star$: 0.5, 0.6 and 0.7.
   (a) The data sets were generated including fixed effects. (b) The matrix $\Z$ used to generate data sets comes from the IMAGEN data. 
The black dots correspond to the mean of $\hat{\eta}$ over 10 replications and the crosses are the real value of $\eta^\star$.}
  \label{fig:fix_vrai}
\end{figure}

We conducted some additional simulations in order to see the impact of the linkage disequilibrium, that is the possible correlations between the columns of $\Z$. Indeed, in the previous numerical study, we generated a matrix $\W$ with independent columns. The matrix $\W$ that we use now to generate the observations is the one from our genetic data set, except that we truncated it in order to have $n=2000$ and $N=100000$. The results of this additional study are presented in Figure \ref{fig:fix_vrai} (b). We can see that they are similar to those obtained previously in Figure 
\ref{fig:seuil_076}, which means that our method does not seem to be sensitive to the presence of correlation between the columns of $\W$.

 \subsubsection{Computational times}

The implementation that we propose in the R package EstHer is very efficient since it only takes 45 seconds for estimating the heritability
and 300 additional seconds to compute the associated 95\% confidence interval. These results have been obtained 
with a computer having the following configuration: RAM 32 GB, CPU 4 $\times$ 2.3 GHz.


\subsubsection{Recovering the support}

\begin{figure}[!ht]
  \centering
\begin{tabular}{cc}
  \includegraphics[width=.45\textwidth] {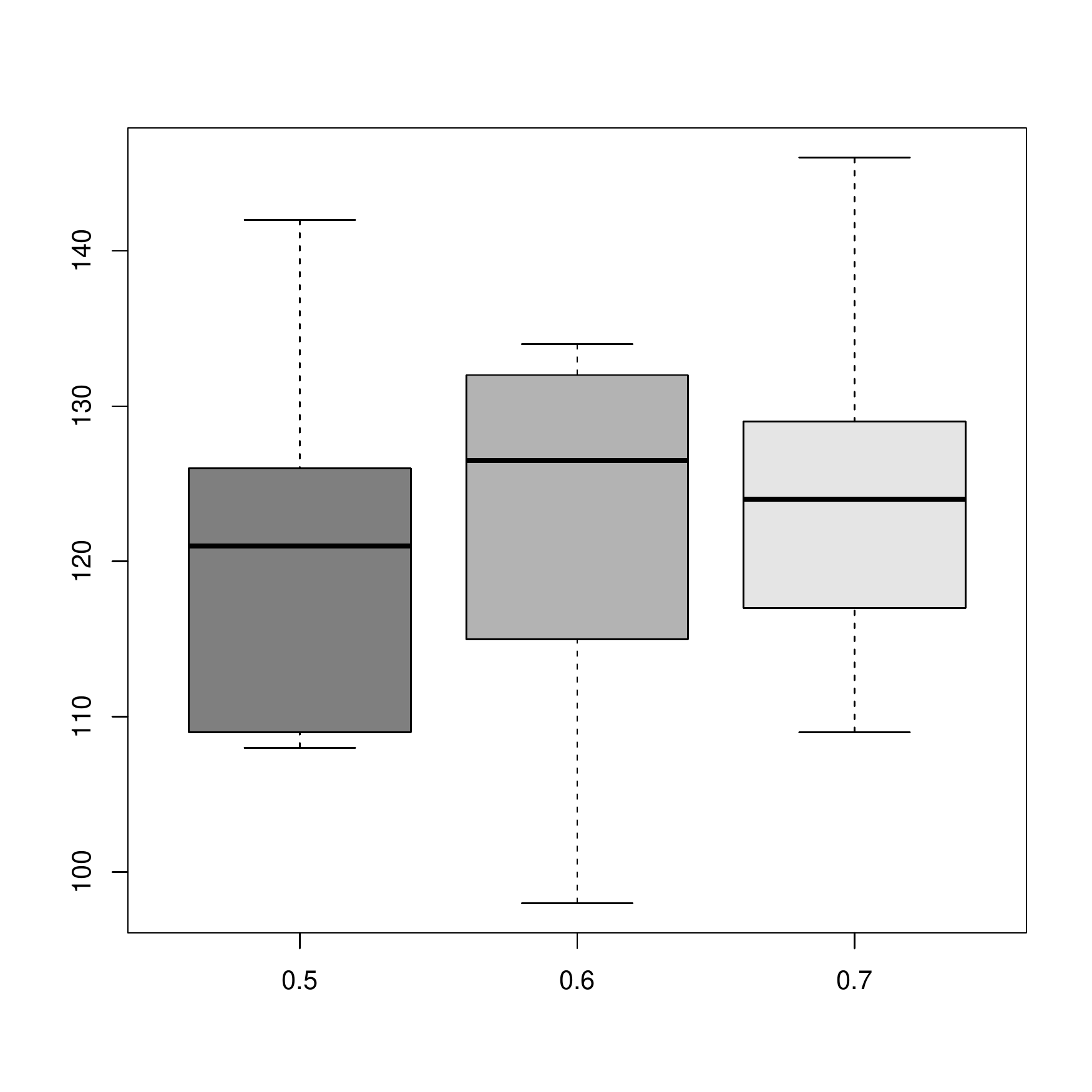}&
  \includegraphics[width=.45\textwidth]{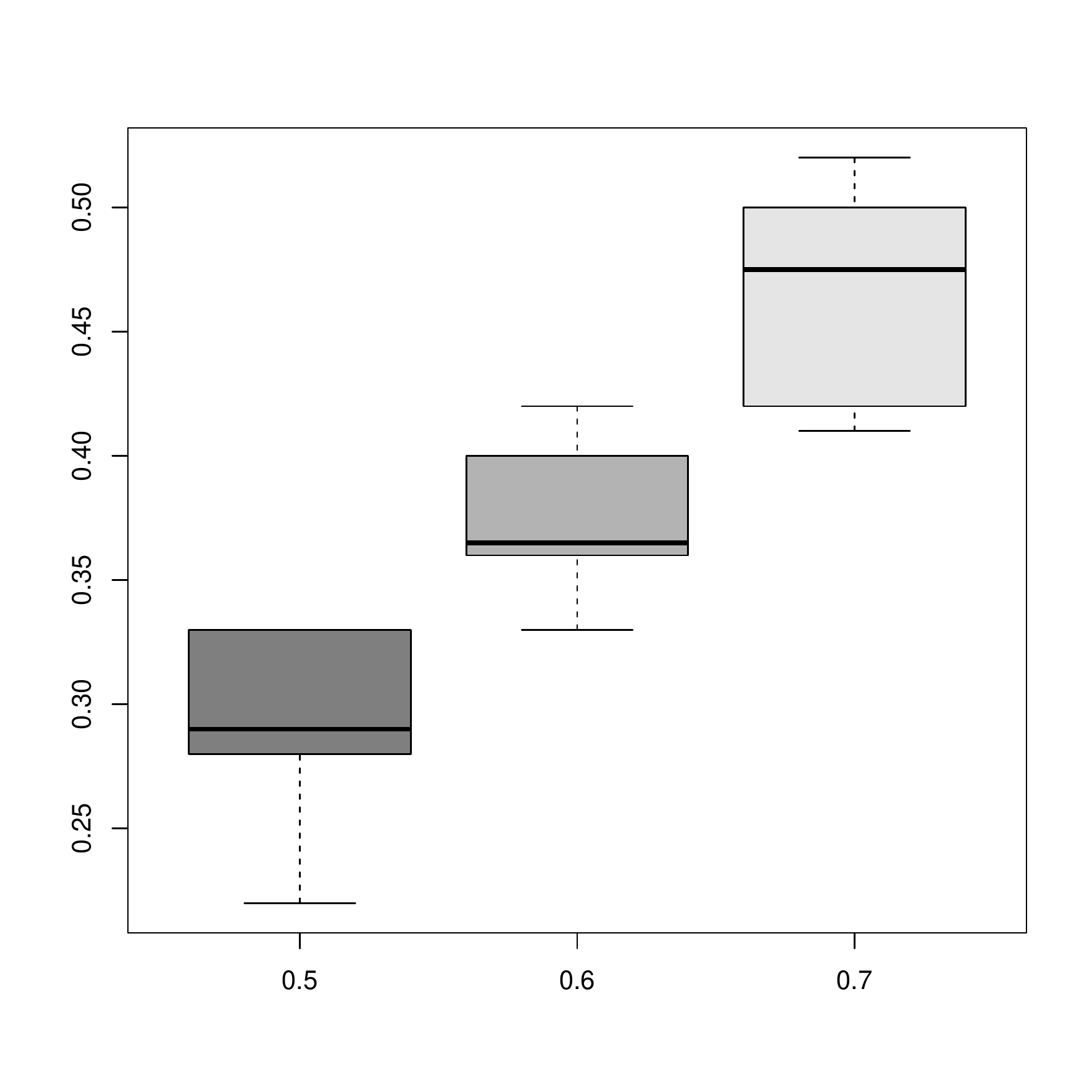}\\
(a) & (b)\\
\end{tabular}
\caption{(a) Boxplots of the length of the set of selected variables with EstHer  for 40 repetitions. The real number of non null components is $100$. (b) Boxplots of the proportion of the real non null components captured in the set of selected variables.}
  \label{fig:long_ind}
\end{figure}

\begin{figure}[!ht]
  \centering
\begin{tabular}{cc}
  \includegraphics[trim=10mm 15mm 0mm 20mm,clip,width=.450\textwidth] {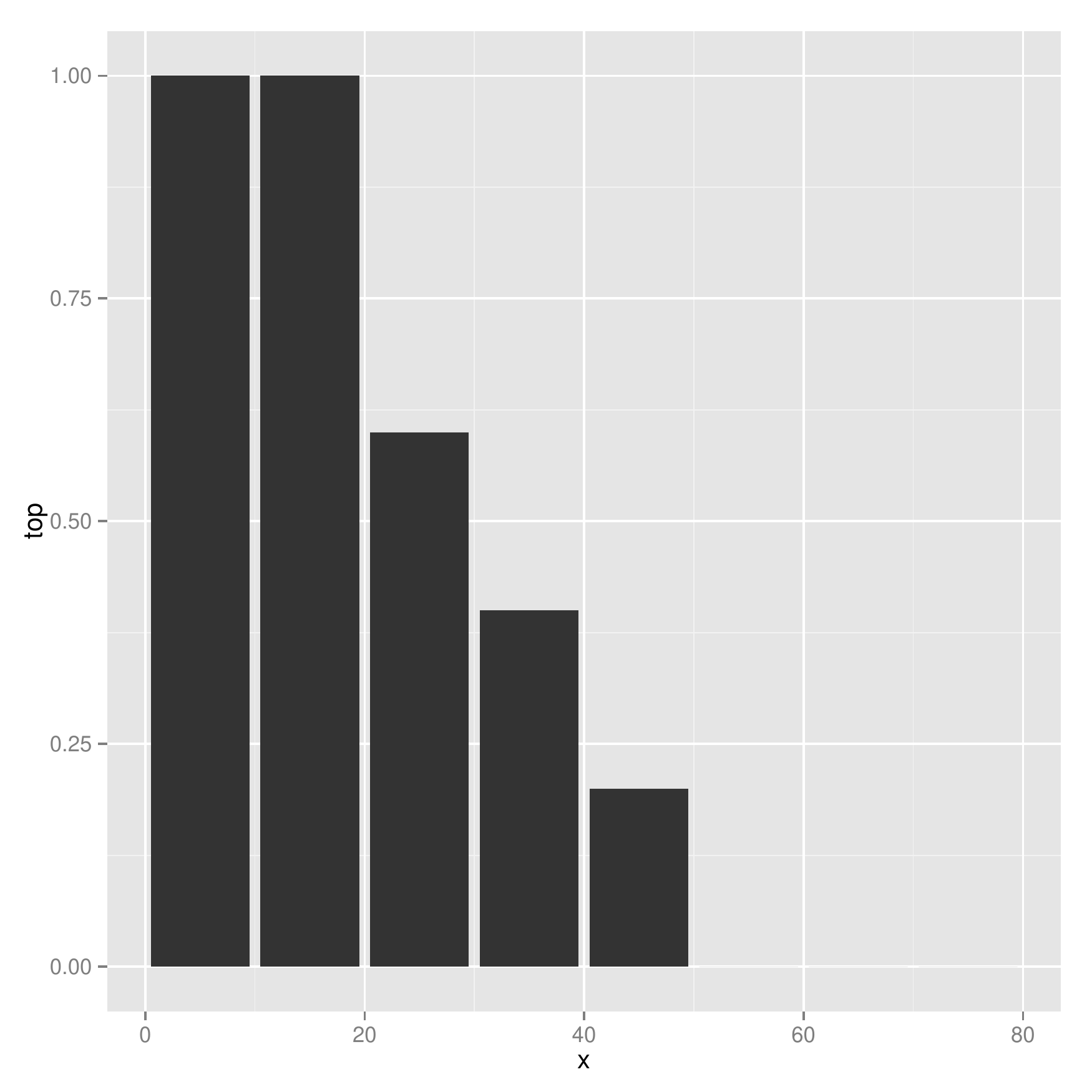}
  &\includegraphics[trim=10mm 15mm 0mm 20mm,clip,width=.45\textwidth]{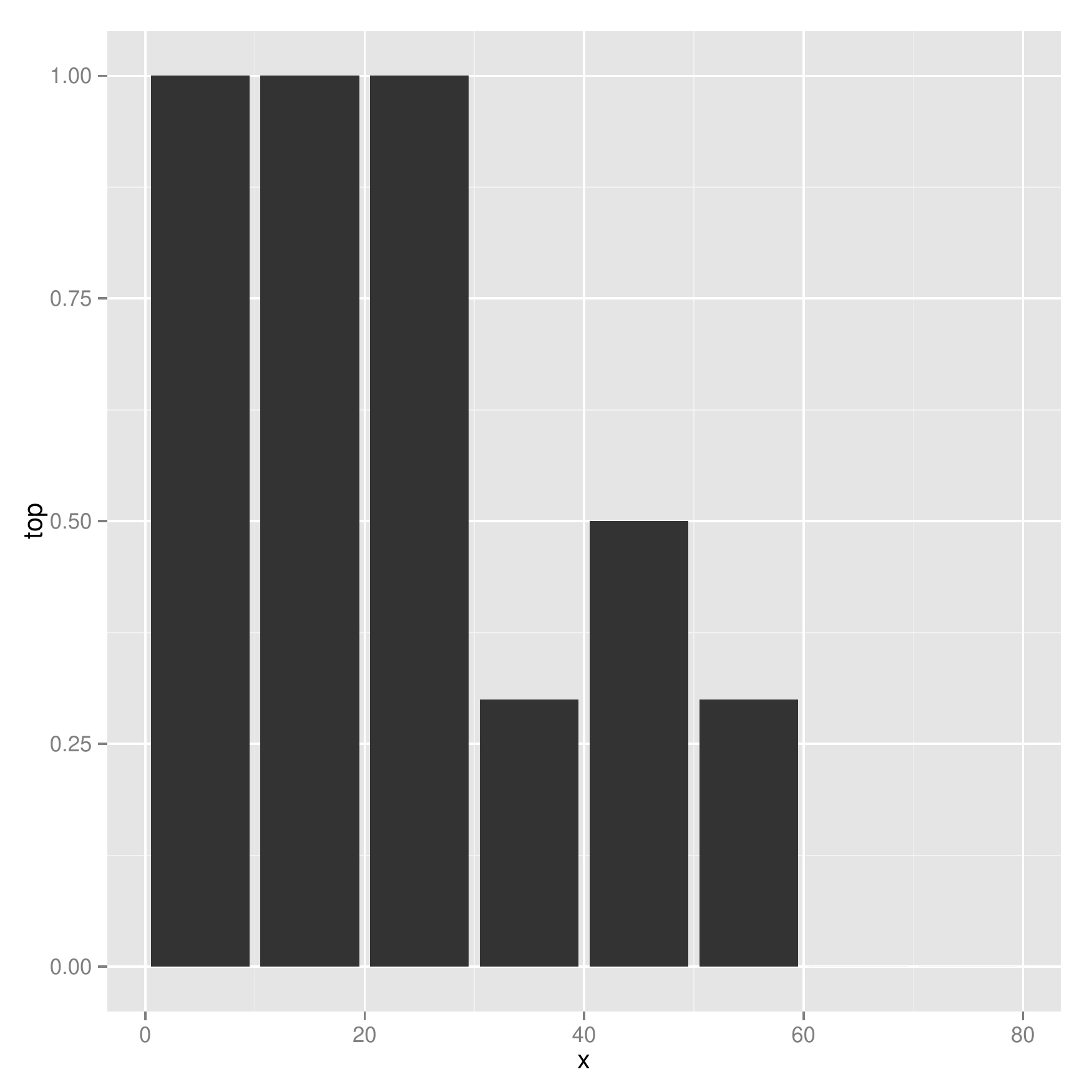} \\
(a) & (b)\\
\end{tabular}
\includegraphics[trim=10mm 15mm 0mm 20mm,clip,width=.45\textwidth]{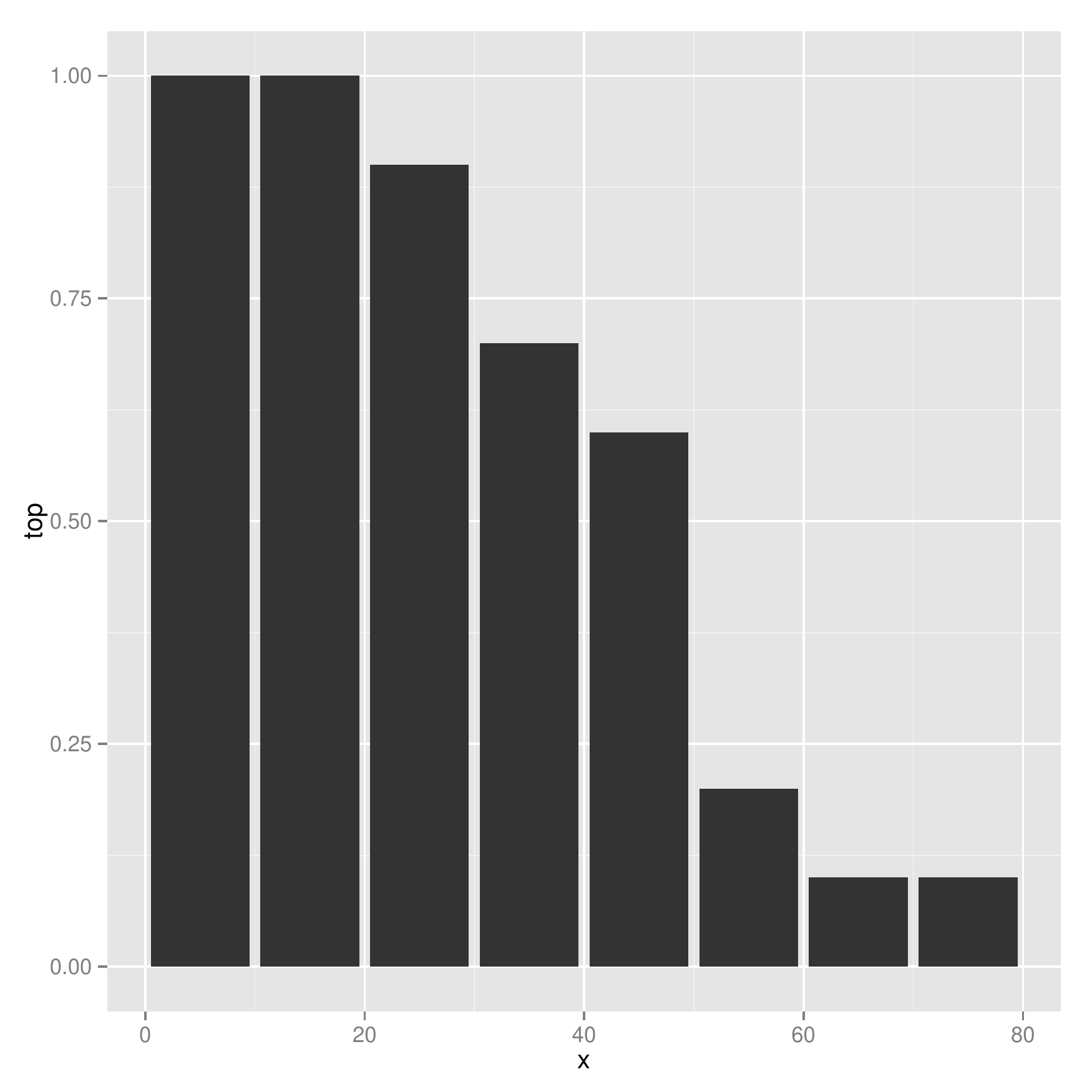}\\
(c)\\
\caption{ Barplots of the proportion of components found by our method as function of the most efficient variables. For example, the first bar is the proportion of the 10 \% higher components that we captured with our selection method. The histograms are displayed for several values of $\eta^{\star}$: 0.5 (a), 0.6 (b), 0.7 (c).}
  \label{fig:histo_grd_comp}
\end{figure}

When the number of causal SNPs is reasonably small, our variable selection method is efficient to estimate the heritability and we wonder if it is reliable as well to recover the support of the random effects.
In Figure \ref{fig:long_ind}, we see the proportion of support estimated by our method when there are $100$ causal SNPs: our method selects around $130$ components. We then focus on the proportion of the real support which has been captured by our method: we see that it may change according to  $\eta^\star$. Indeed, the higher  $\eta^\star$, the higher this proportion.
 Nevertheless, even in the worst case, that is  $\eta^\star=0.5$, Figure \ref{fig:histo_grd_comp} shows that even if we keep only $30$\% of the real non null components, we select the most active ones. 

The ability of recovering the support in linear models has been studied by \cite{verzelen:2012} in ultra high dimensional cases. The author shows that with a non null probability, the support cannot be estimated under some numerical conditions on the parameters $q$, $N$ and $n$ (namely if there are considerably more variables $N$ than observations $n$, and if the number of non null components $qN$ is relatively high). In this simulation study, even when we consider small values of $q$ (for instance $q=10^{-3}$, that is $100$ causal SNPs), we are not far from to the ultra high dimensional framework described  in \cite{verzelen:2012}, which can explain the difficulties to recover the full support. 

\subsection{\textcolor{black}{Results when the number of causal SNPs is high}}

\textcolor{black}{In subsection \ref{subsec:res_q_small} we show the performance of our method in the case where the proportion of causal SNPs $q$ is small, that is around $10^{-3}$. In this subsection, we focus on a more polygenic scenario, that includes the cases where thousands of SNPs or ten of thousands of SNPs are causal.}

\subsubsection{\textcolor{black}{Results when there are SNPs with moderate and weak effects}}

\begin{figure}[!ht]
  \centering
\begin{tabular}{cc}
$\eta^\star=0.4$ & $\eta^\star=0.6$ \\
  \includegraphics[trim=10mm 15mm 0mm 20mm,clip,width=.4\textwidth] {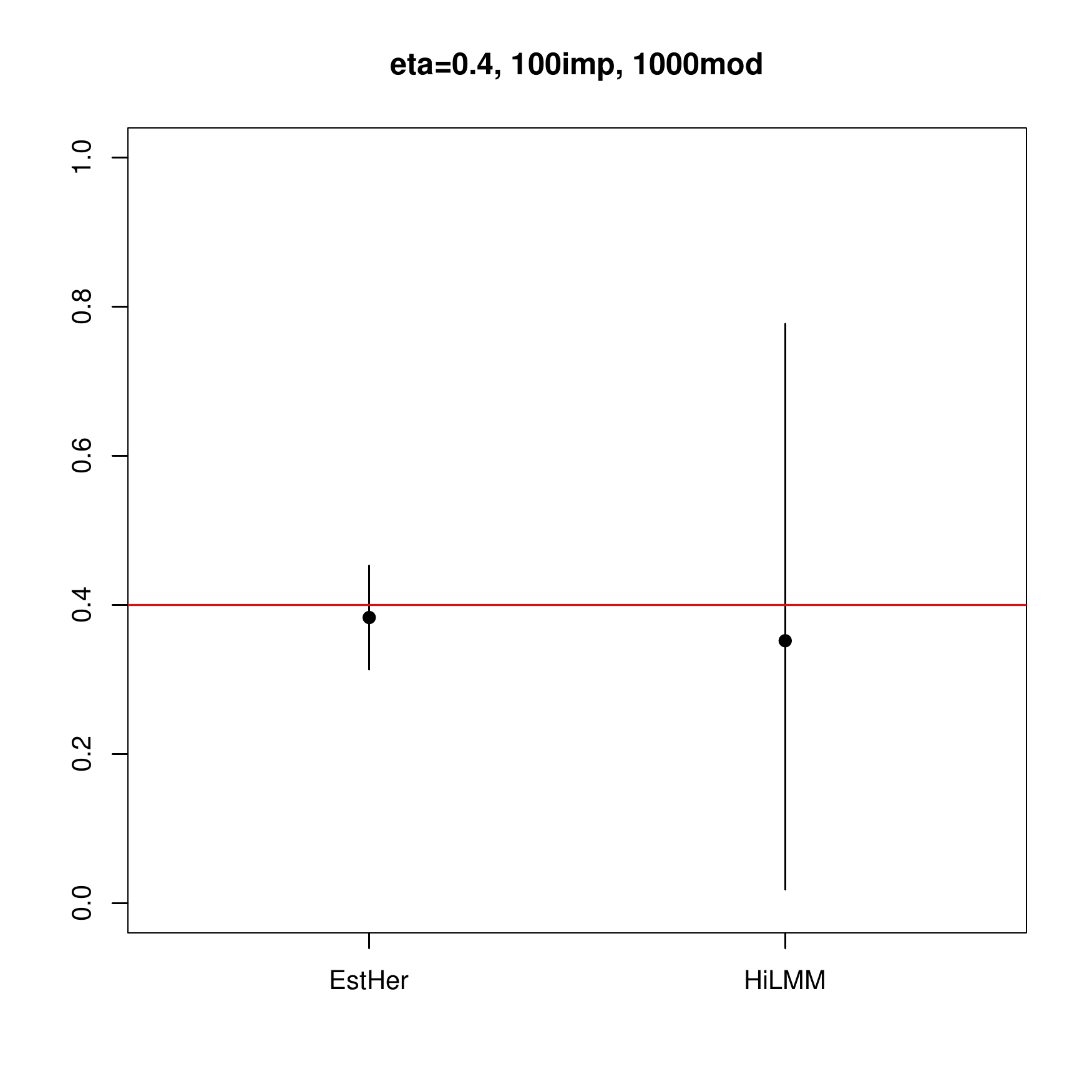}
&  \includegraphics[trim=10mm 15mm 0mm 20mm,clip,width=.4\textwidth]{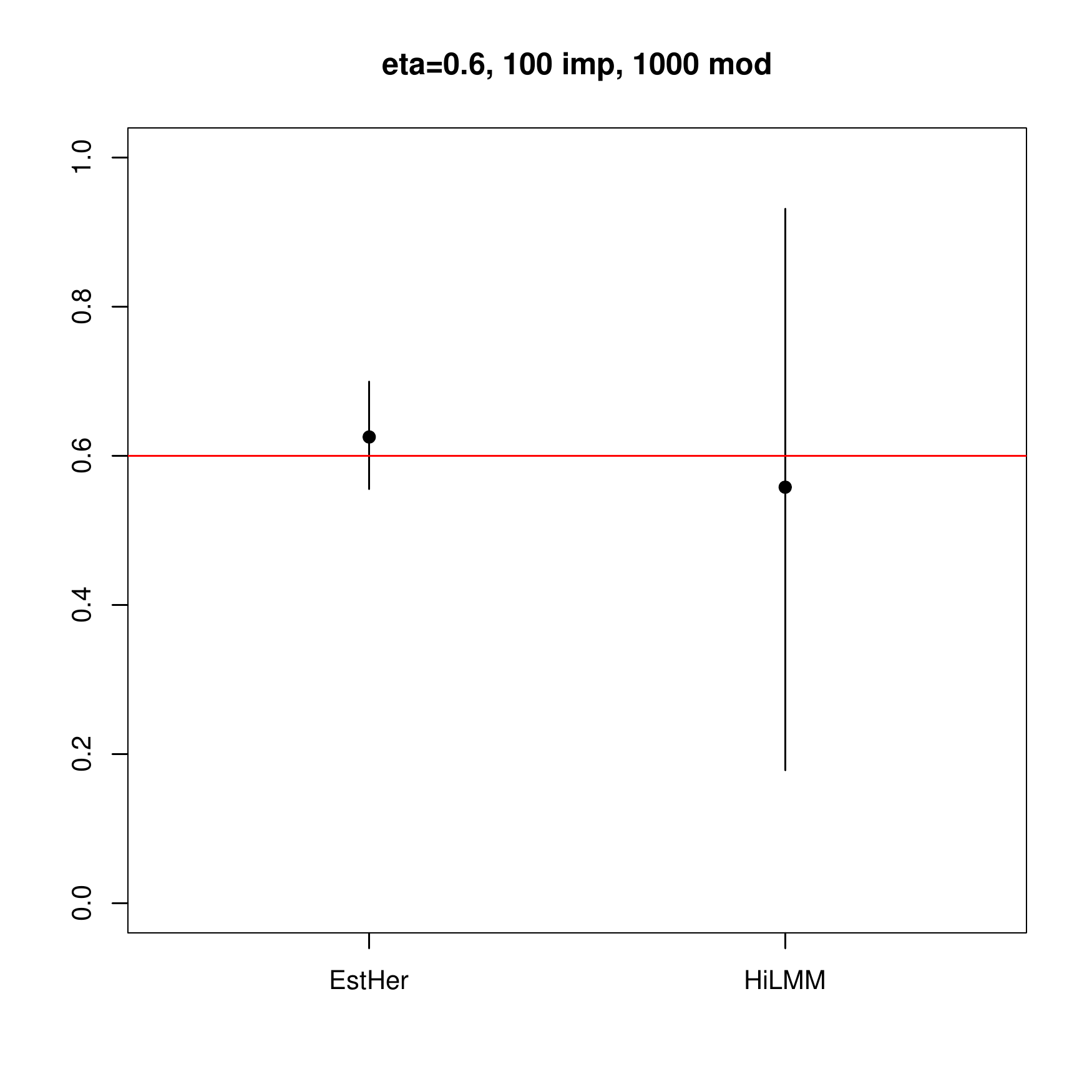} \\
  \includegraphics[trim=10mm 15mm 0mm 20mm,clip,width=.4\textwidth] {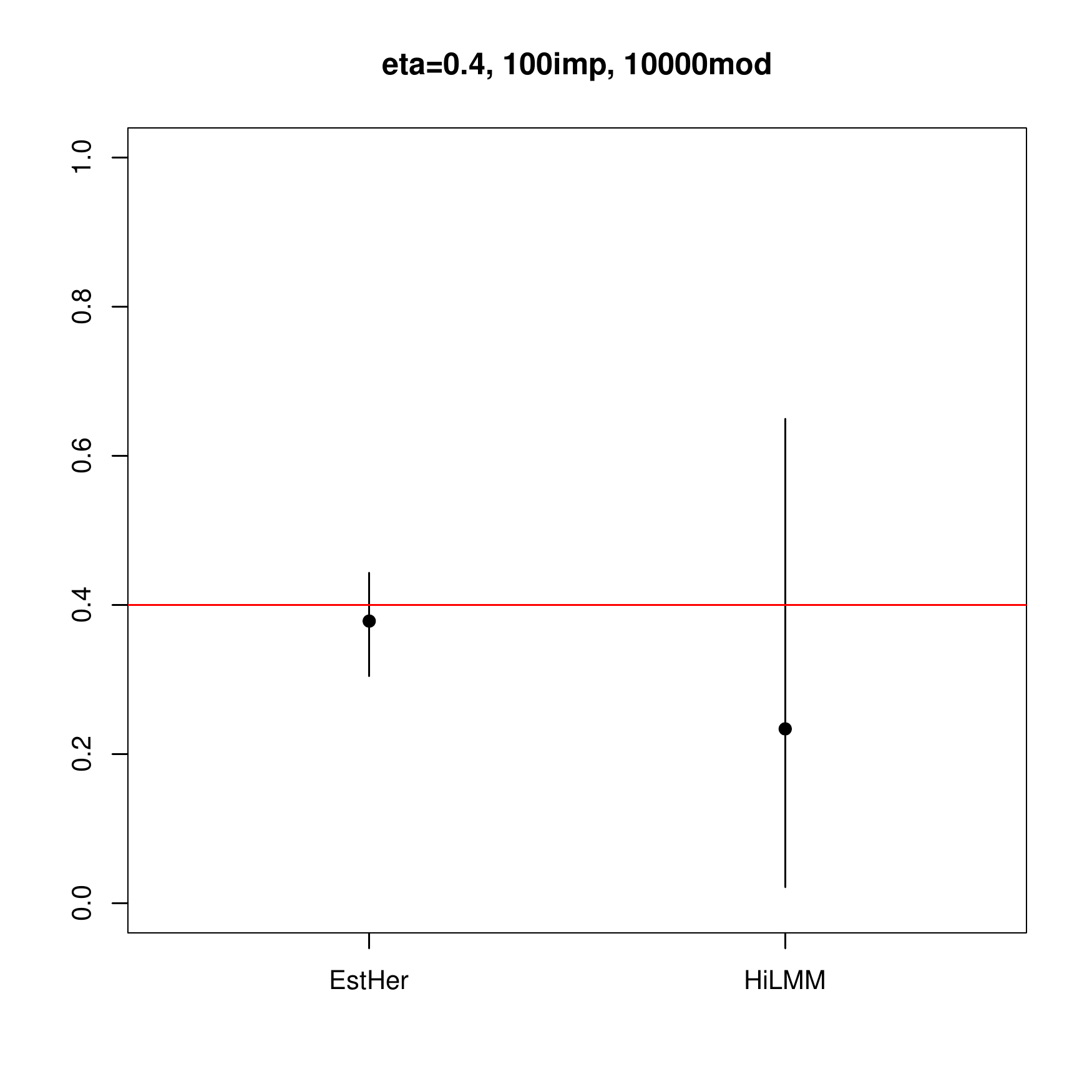}
&  \includegraphics[trim=10mm 15mm 0mm 20mm,clip,width=.4\textwidth]{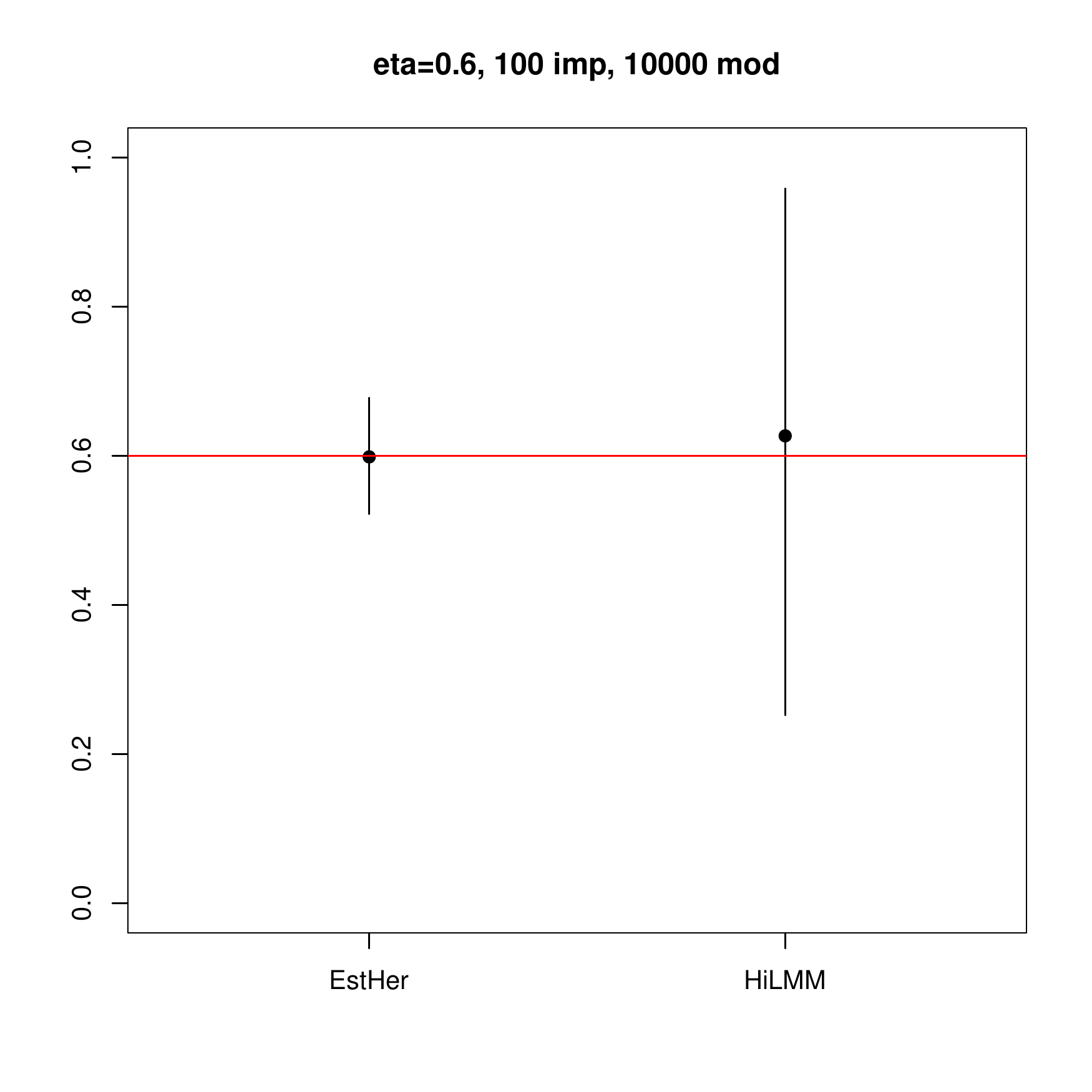}
\\
\end{tabular}
\caption{Results of HiLMM and EstHer when there are a few causal SNPs with moderate effects and a lot of SNPs with small effects. The proportion of each is 100 out of 1000 (up) and 100 out of 10000 (bottom), with $\eta^\star=0.4$ and $0.6$.}
\label{fig:EstHer_mix}
\end{figure}

\textcolor{black}{We first focus on the statistical performance of EstHer when there are a lot of SNPs (1000 or 10000) with small effects (for example, that explain 5\% of the phenotypic variations), and a small number (around 100) with moderate effects. We can see from
Figure \ref{fig:EstHer_mix} that, in this case, EstHer provides unbiased estimations with a small variance.}

\subsubsection{\textcolor{black}{Results when all SNPs have moderate effects}}

\begin{figure}[!ht]
  \centering
\begin{tabular}{cc}
$\eta^\star=0.4$ & $\eta^\star=0.6$ \\
  \includegraphics[trim=10mm 15mm 0mm 20mm,clip,width=.4\textwidth] {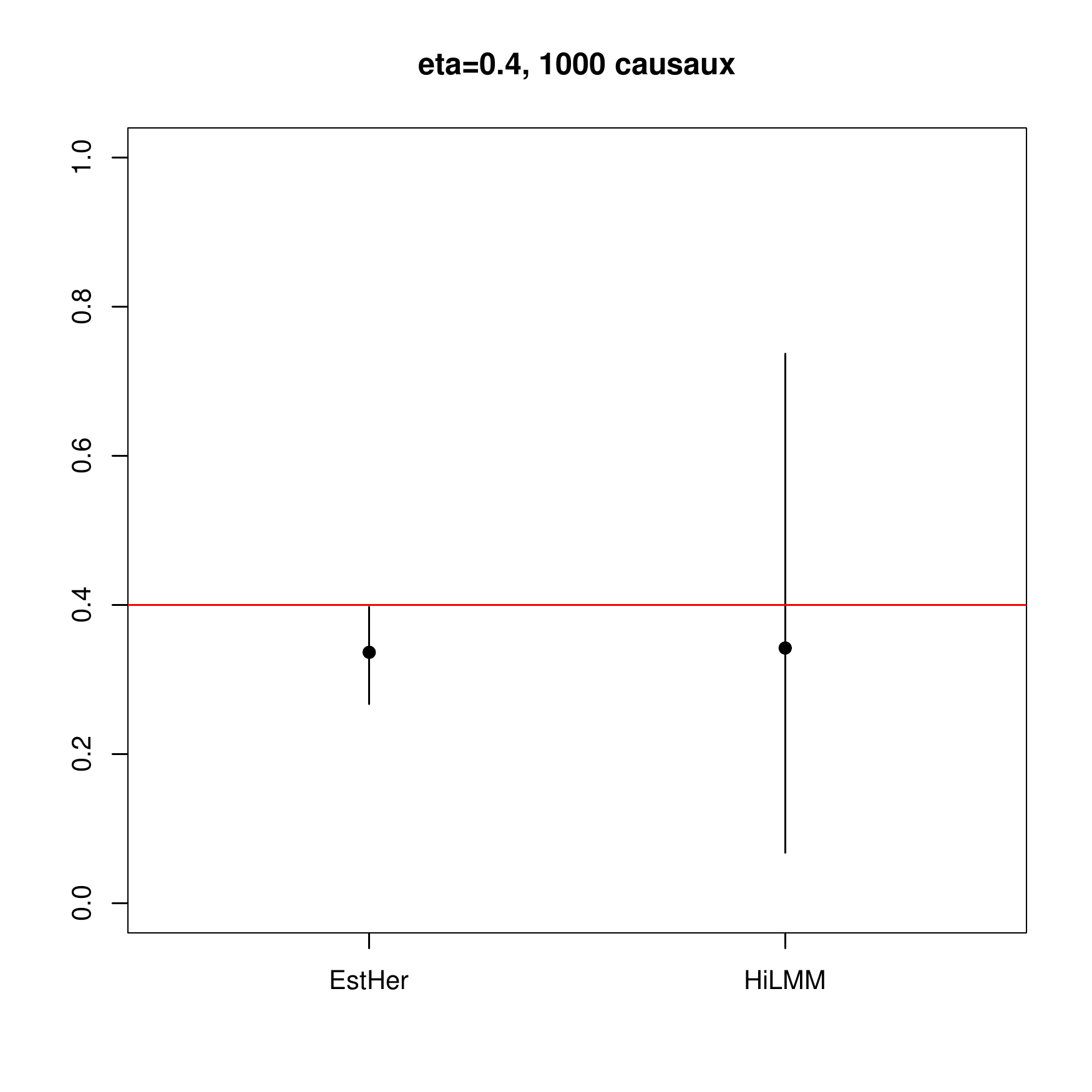}
&  \includegraphics[trim=10mm 15mm 0mm 20mm,clip,width=.4\textwidth]{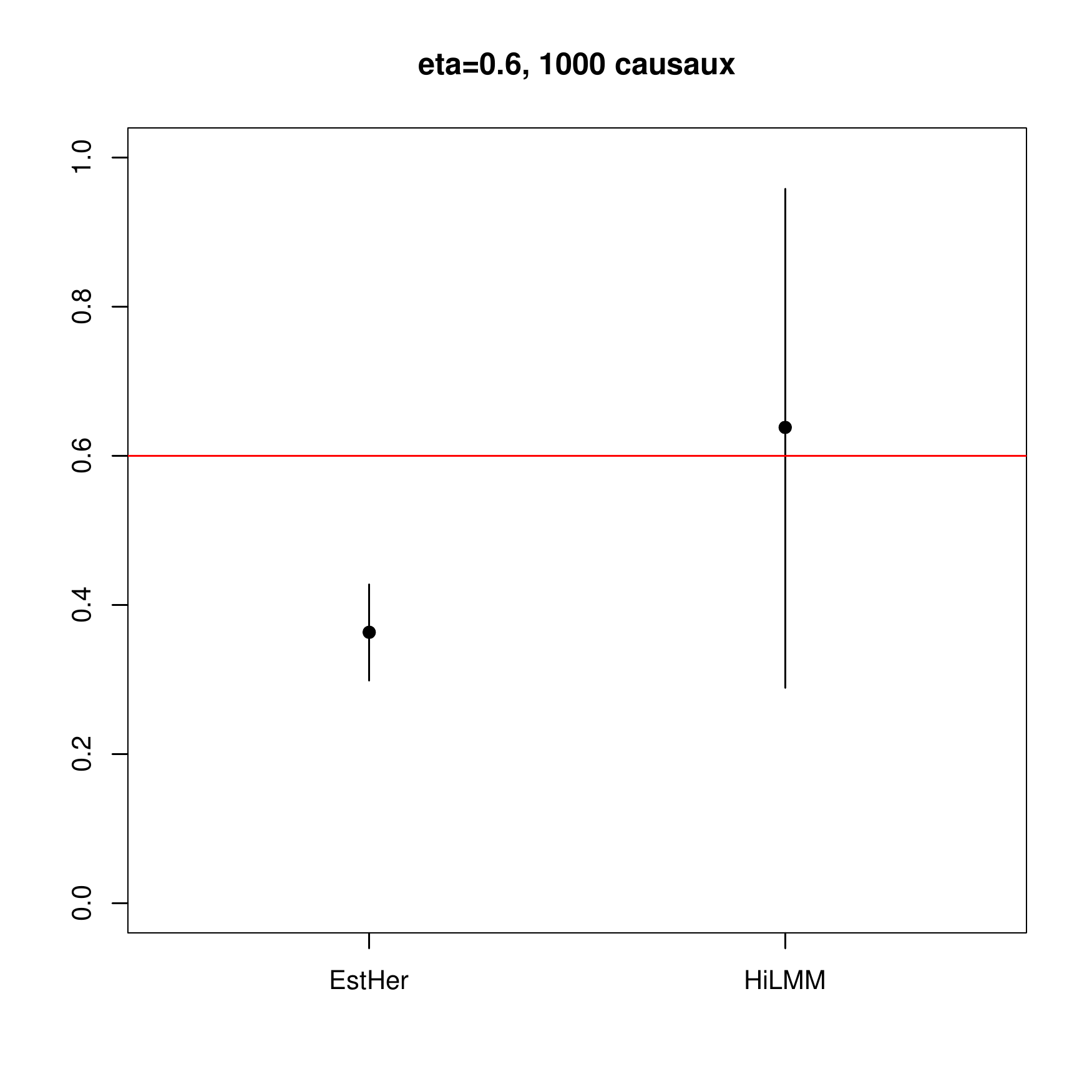} \\
  \includegraphics[trim=10mm 15mm 0mm 20mm,clip,width=.4\textwidth] {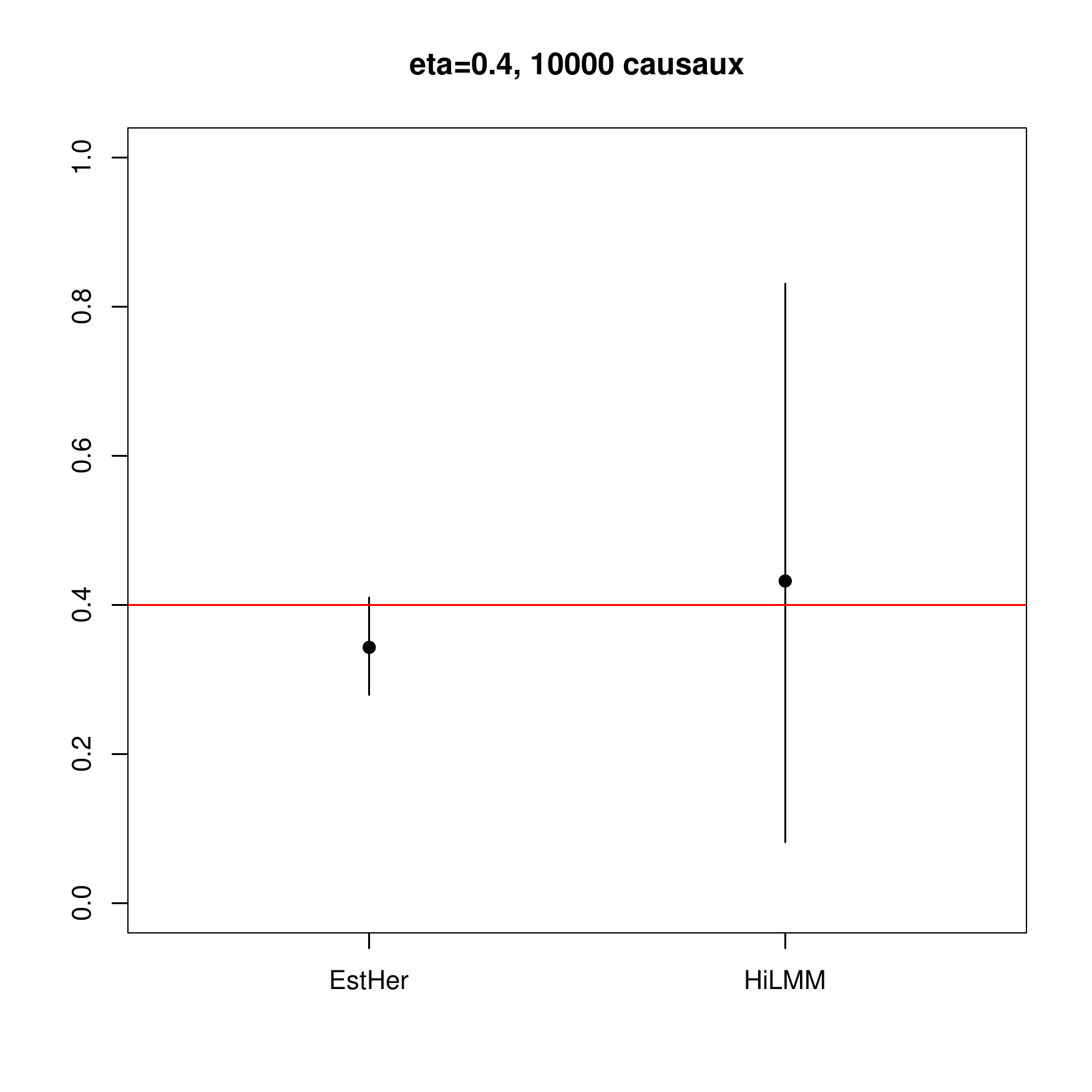}
&  \includegraphics[trim=10mm 15mm 0mm 20mm,clip,width=.4\textwidth]{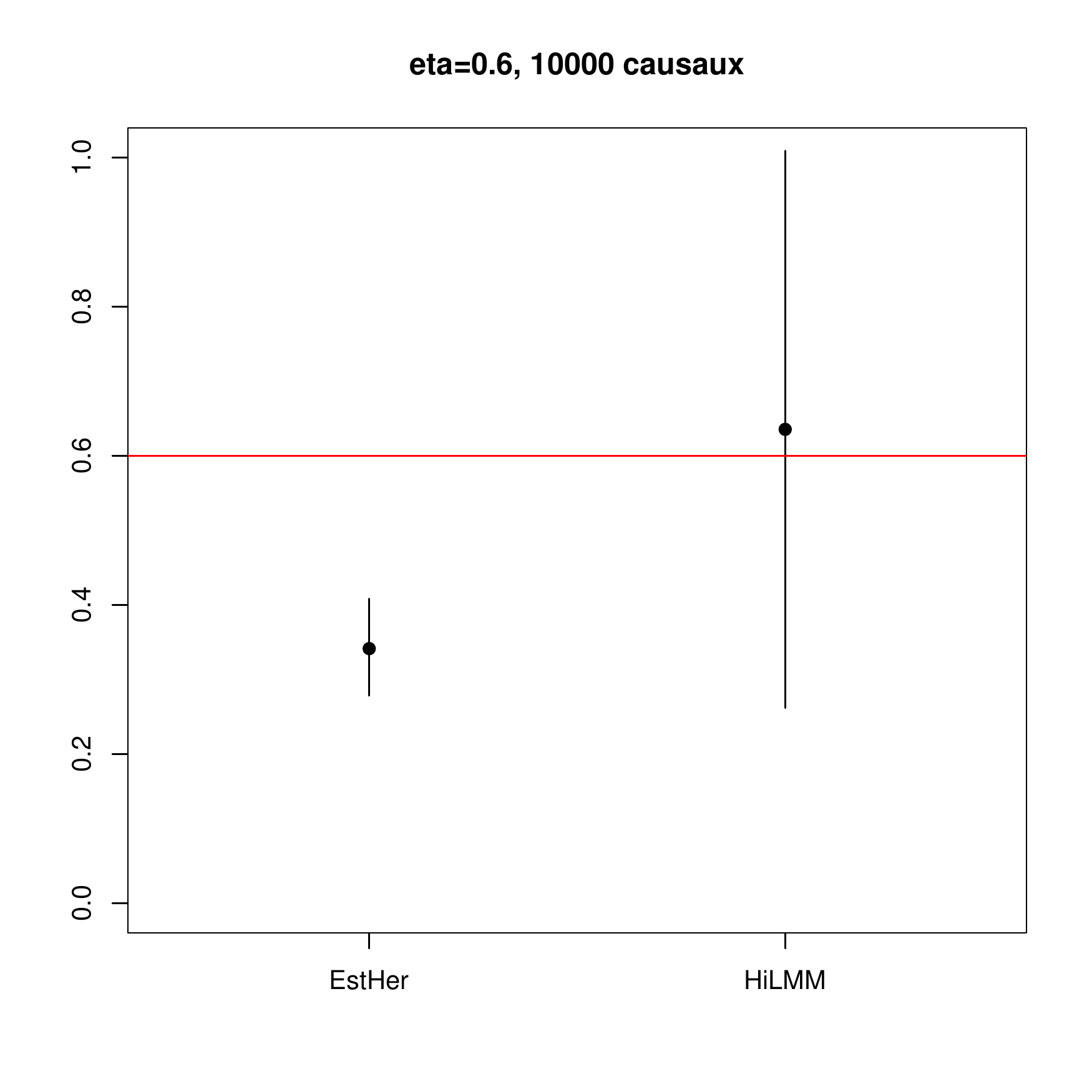}
\\
\end{tabular}
\caption{Results of HiLMM and EstHer for 1000 (up) and 10000 (bottom) causal SNPs and for $\eta^\star=0.4$ and $0.6$.}
\label{fig:EstHer_fails}
\end{figure}

\begin{figure}[!ht]
\begin{tabular}{cc}
  \includegraphics[trim=10mm 15mm 0mm 20mm,clip,width=.4\textwidth] {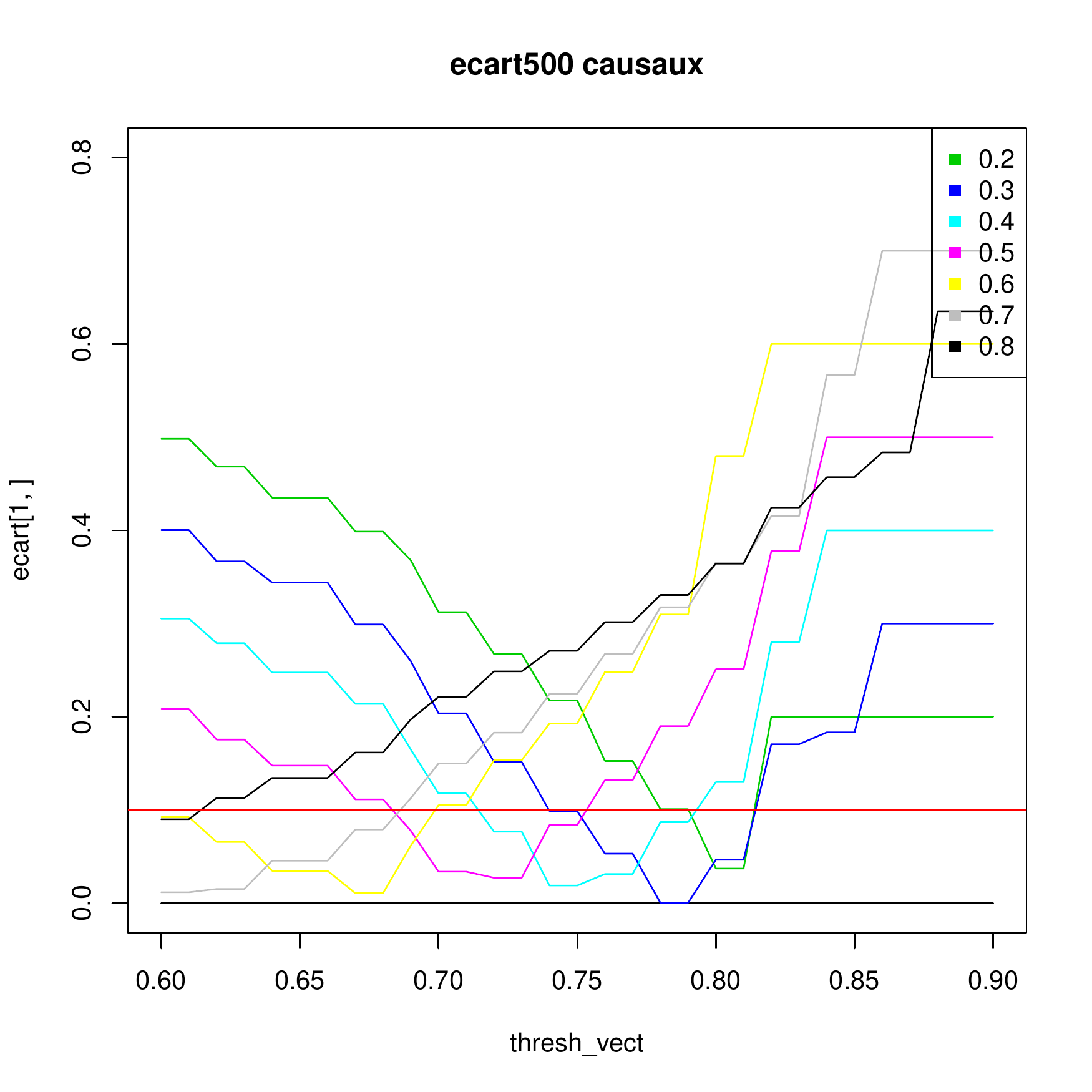}
&  \includegraphics[trim=10mm 15mm 0mm 20mm,clip,width=.4\textwidth]{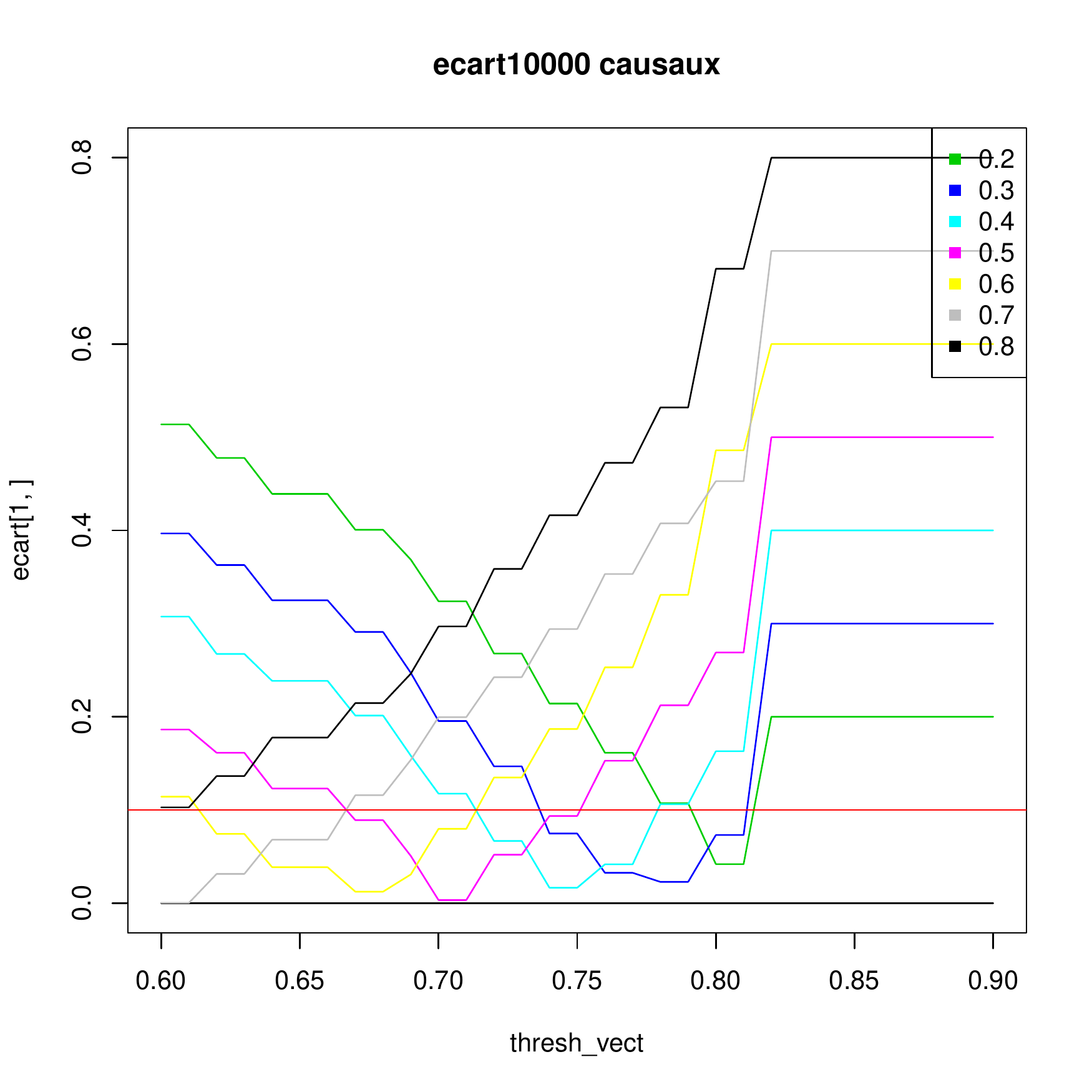} \\
1000 causal SNPs & 10000 causal SNPs\\
\end{tabular}
\caption{\textcolor{black}{Absolute difference} $|\eta^\star-\hat{\eta}|$ for thresholds from 0.6 to 0.9 and for 1000 (left) and 10000 (right) causal SNPs.}
\label{fig:no_common_threshold}
\end{figure}

\textcolor{black}{If all causal SNPs have moderate effects and if the number of these causal SNPs is high, namely greater than 1000, EstHer underestimates the heritability. These results are displayed in Figure \ref{fig:EstHer_fails}. Moreover, we can see from Figure 
\ref{fig:no_common_threshold} that there is no threshold choice that can provide accurate estimations of heritability for all values 
of $\eta^\star$.}

\section{\textcolor{black}{A criterion to decide whether we should apply EstHer or HiLMM}}
\label{sec:criterion}

\begin{figure}[!ht]
  \centering
\begin{tabular}{ccc}
 $\eta^\star=0.4$ & $\eta^\star=0.5$ & $\eta^\star=0.6$ \\
 \includegraphics[trim=10mm 15mm 0mm 20mm,clip,width=.3\textwidth] {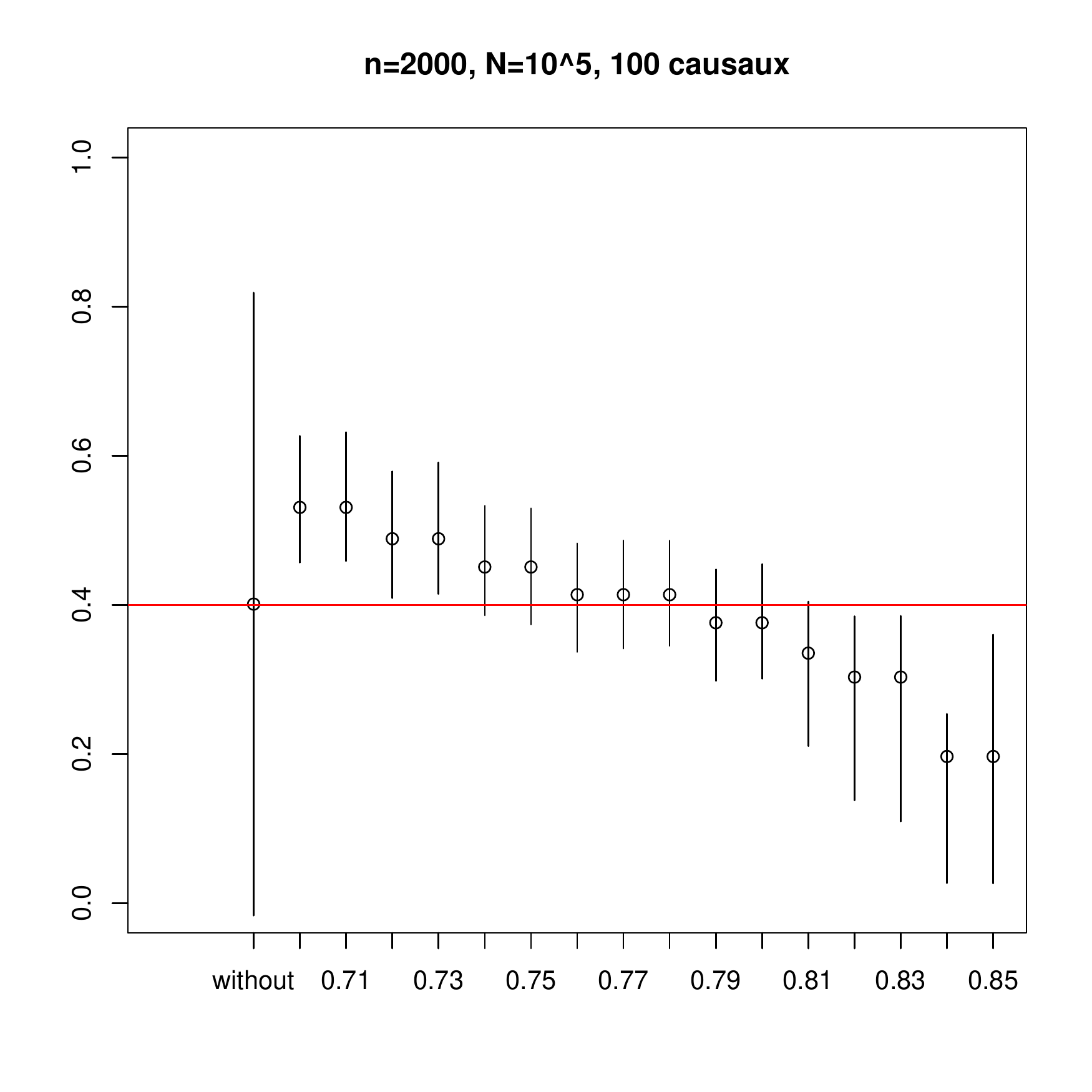}
  &\includegraphics[trim=10mm 15mm 0mm 20mm,clip,width=.3\textwidth]{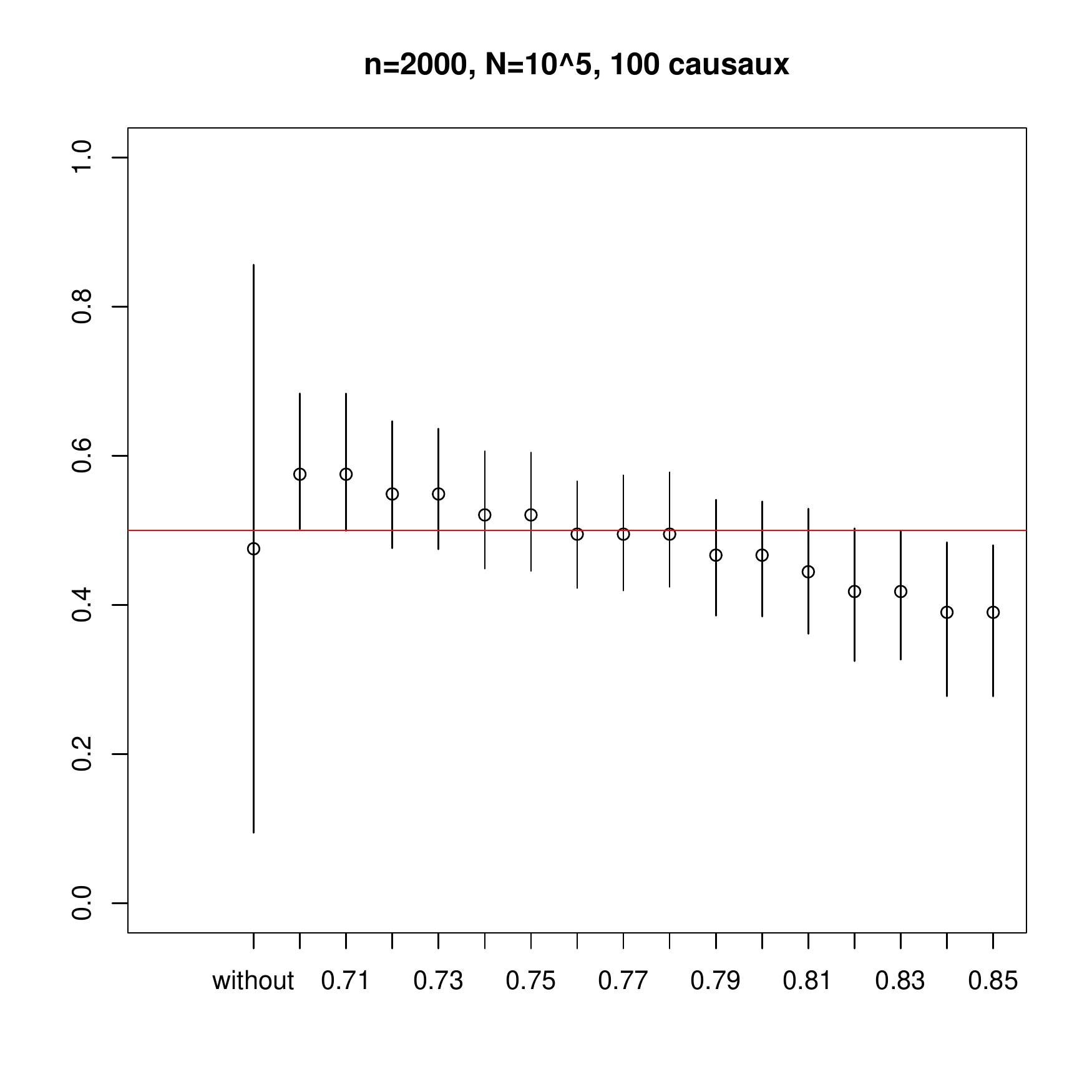} 
&\includegraphics[trim=10mm 15mm 0mm 20mm,clip,width=.3\textwidth]{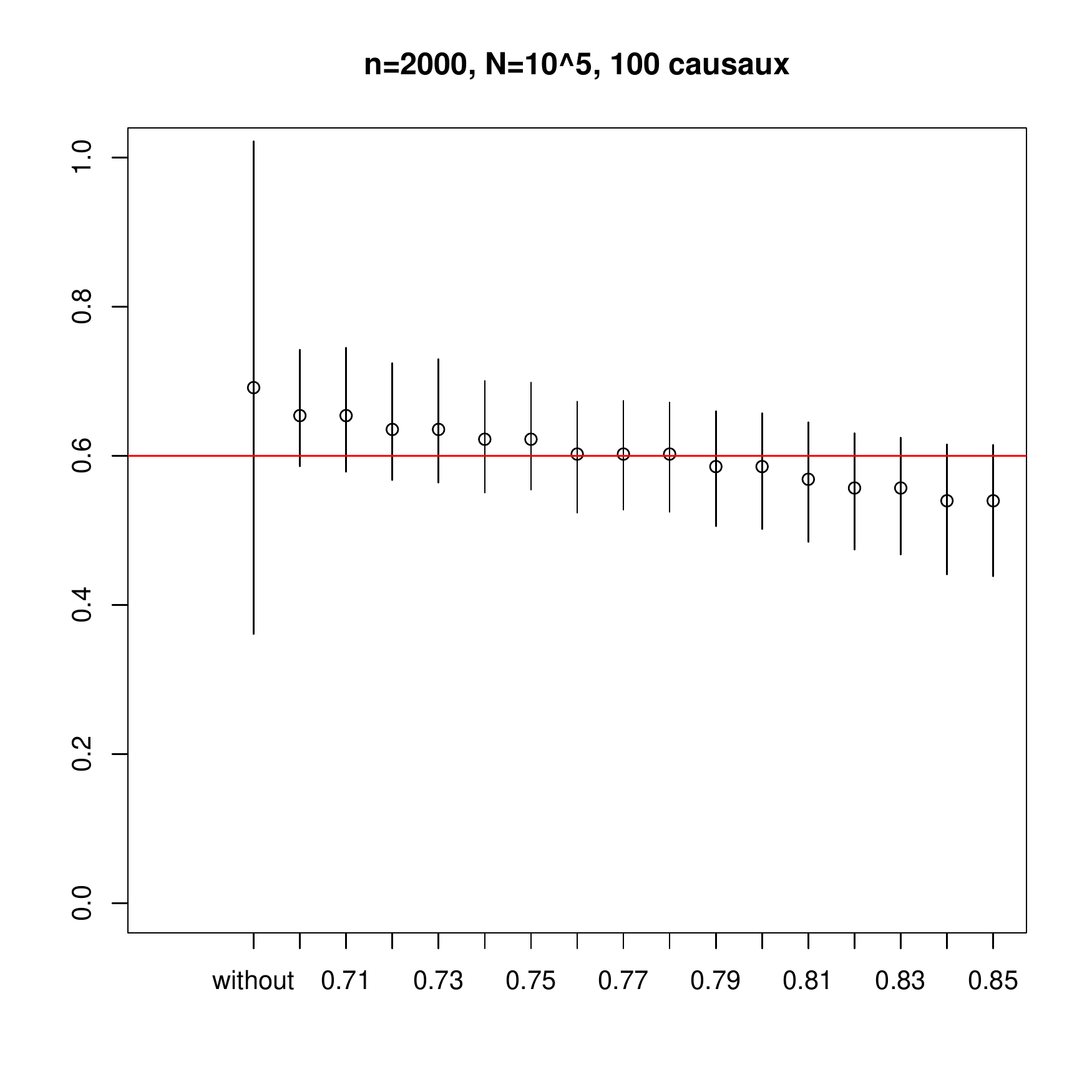}   \\
   \includegraphics[trim=10mm 15mm 0mm 20mm,clip,width=.3\textwidth] {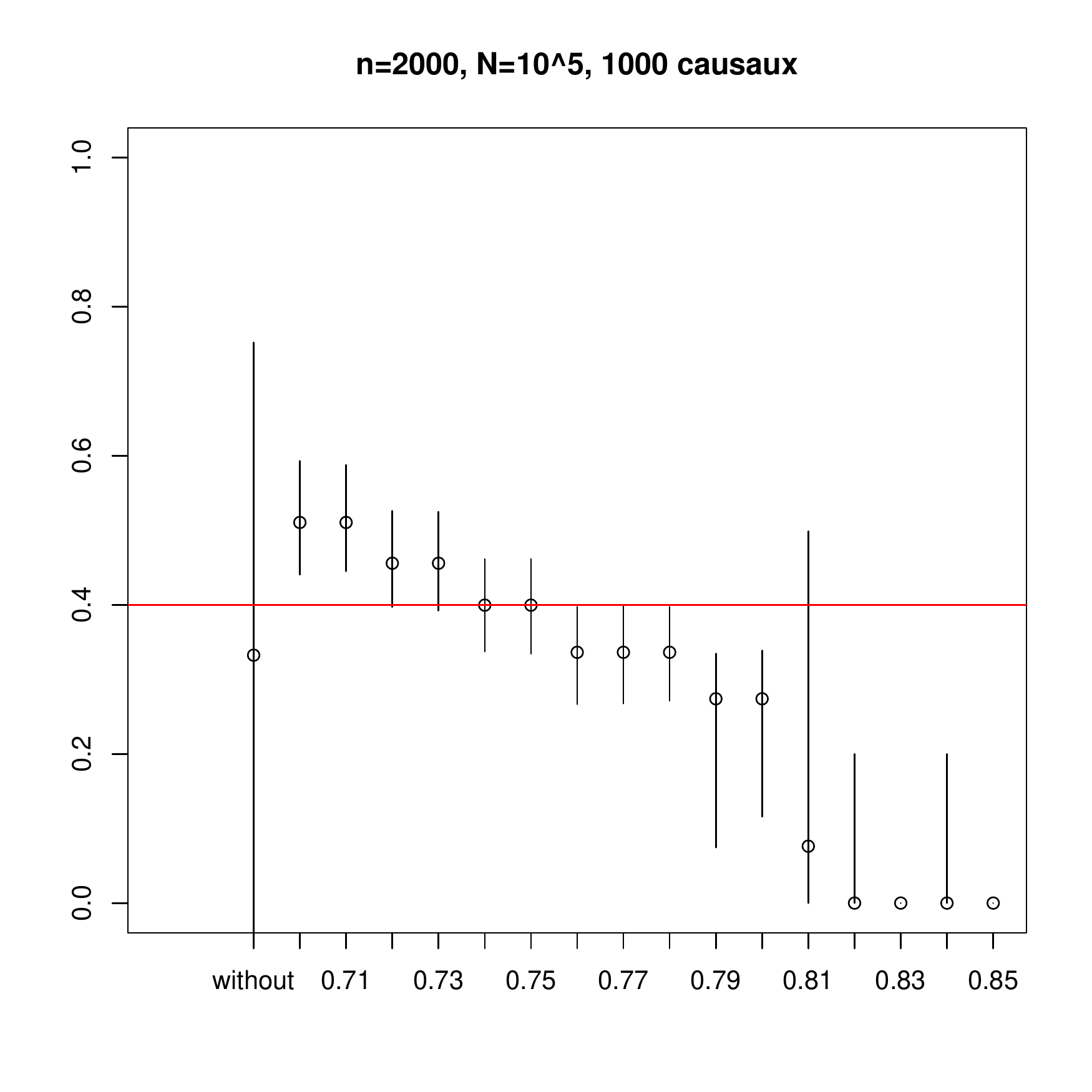}
  &\includegraphics[trim=10mm 15mm 0mm 20mm,clip,width=.3\textwidth]{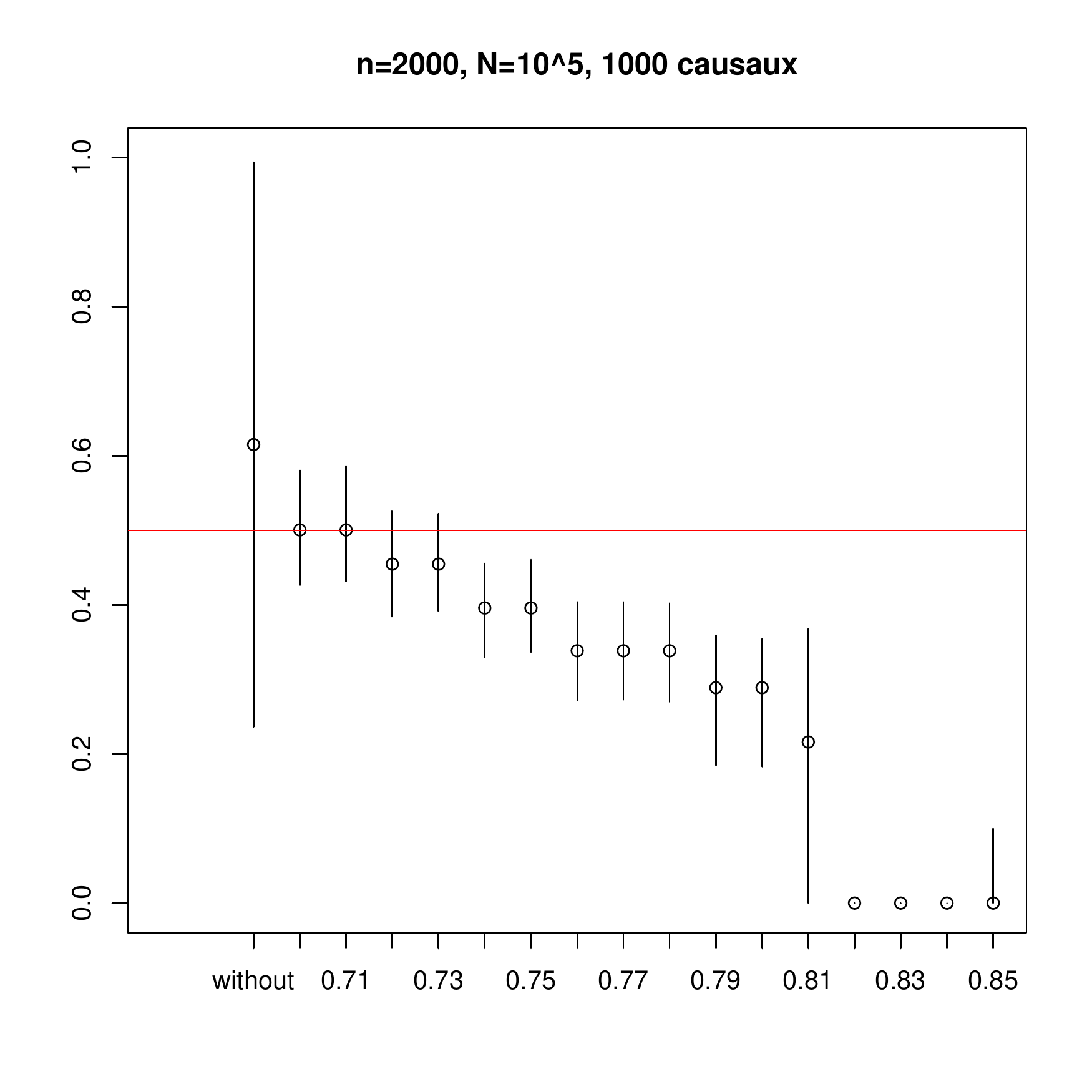} 
&\includegraphics[trim=10mm 15mm 0mm 20mm,clip,width=.3\textwidth]{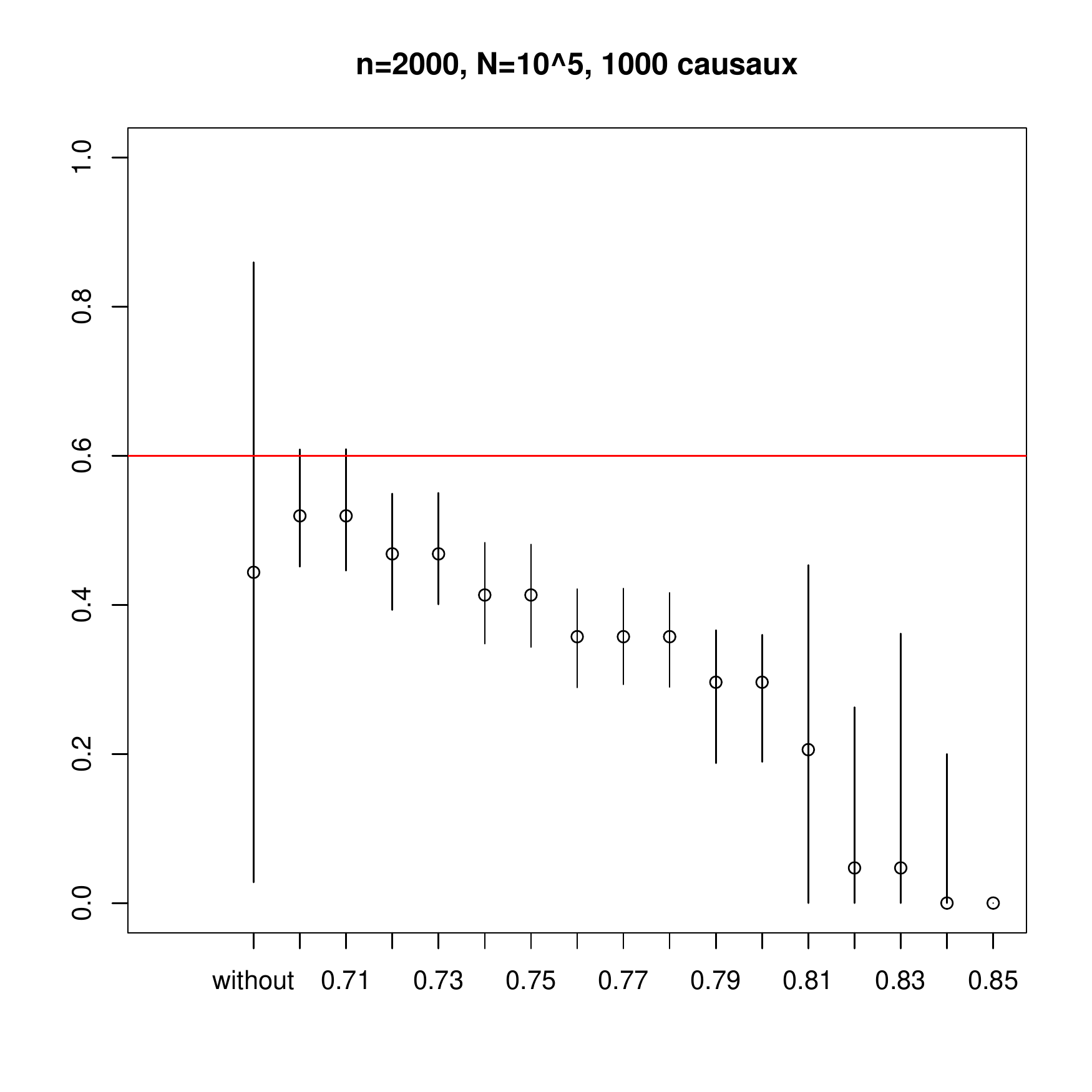}\\
 \includegraphics[trim=10mm 15mm 0mm 20mm,clip,width=.3\textwidth]{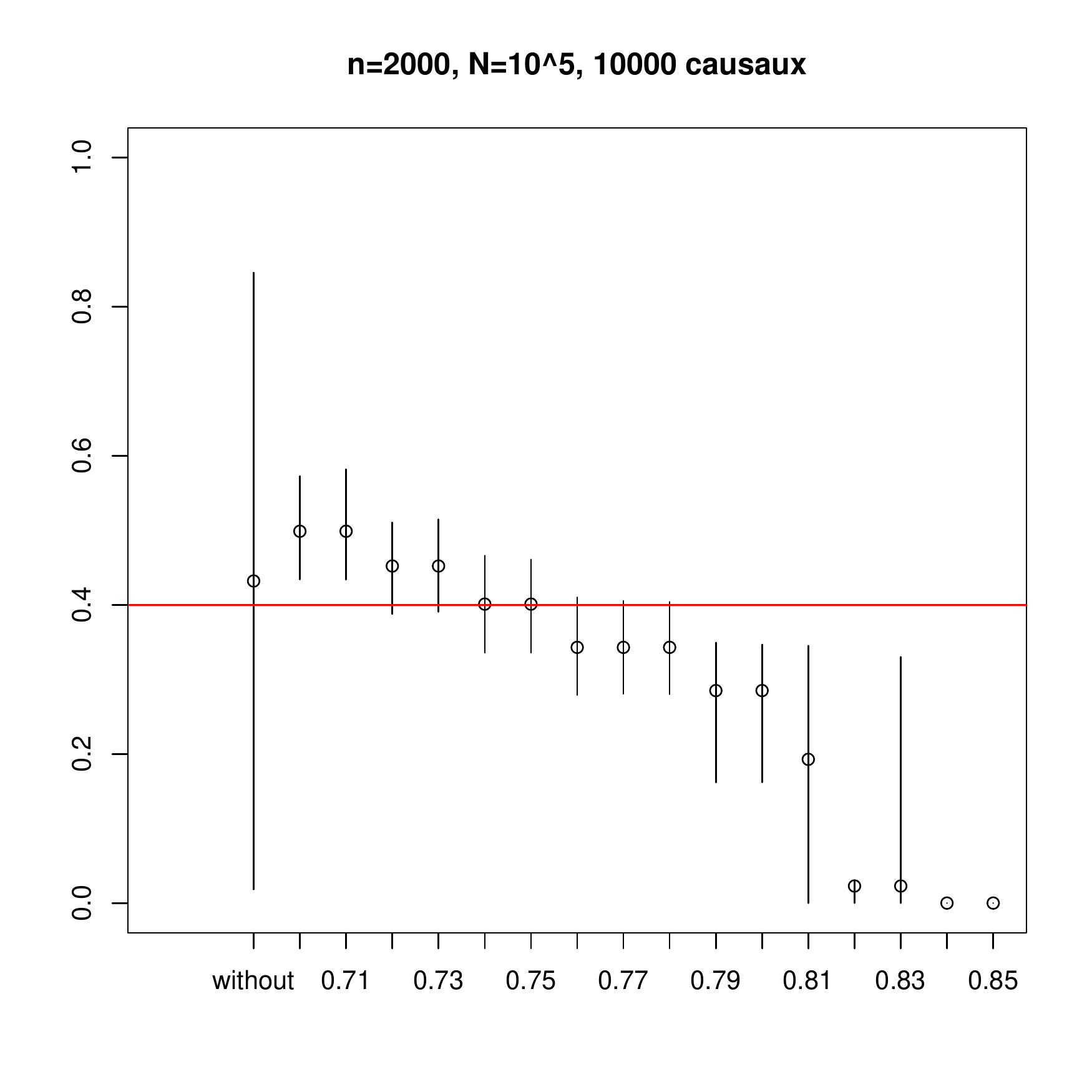}
&\includegraphics[trim=10mm 15mm 0mm 20mm,clip,width=.3\textwidth]{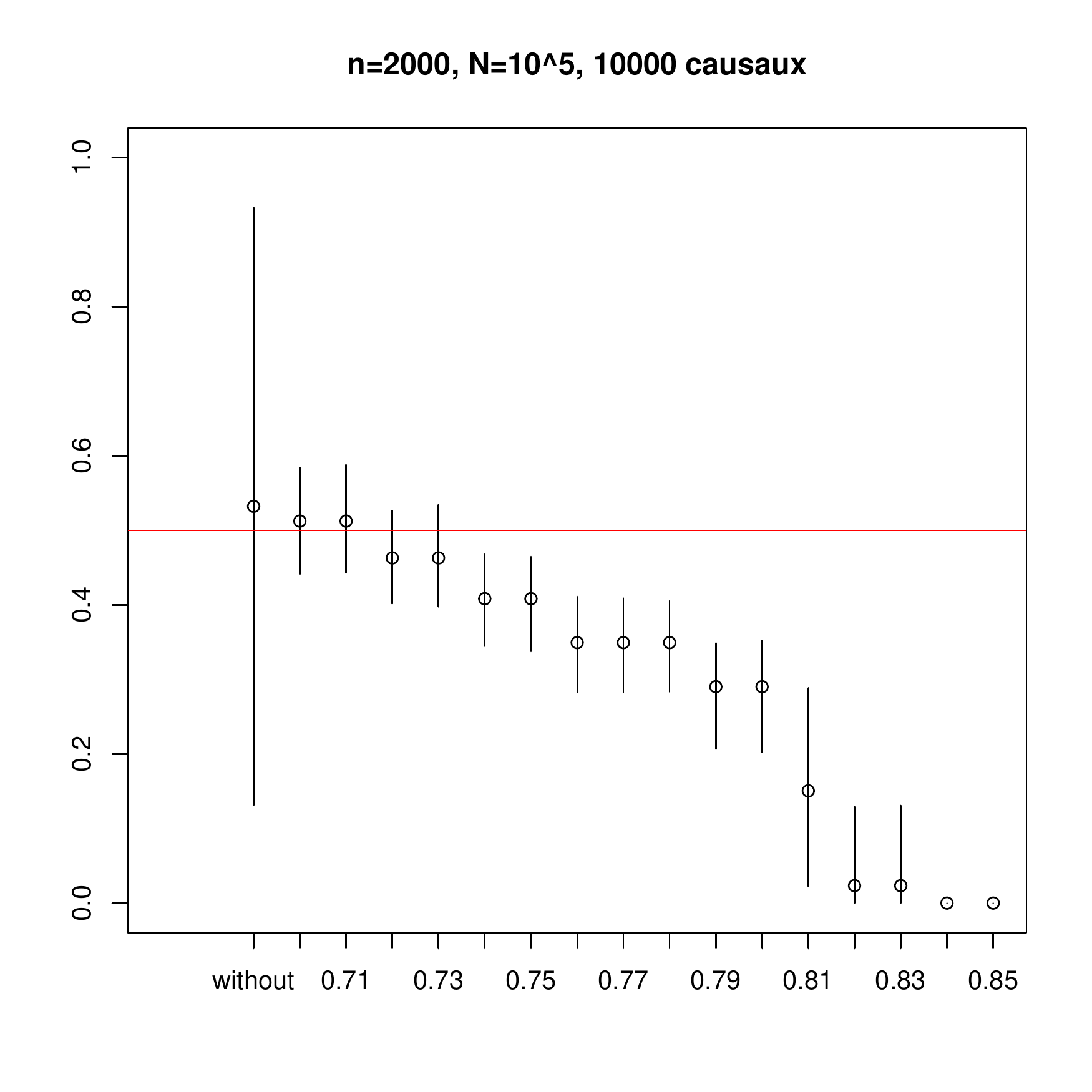}
&\includegraphics[trim=10mm 15mm 0mm 20mm,clip,width=.3\textwidth]{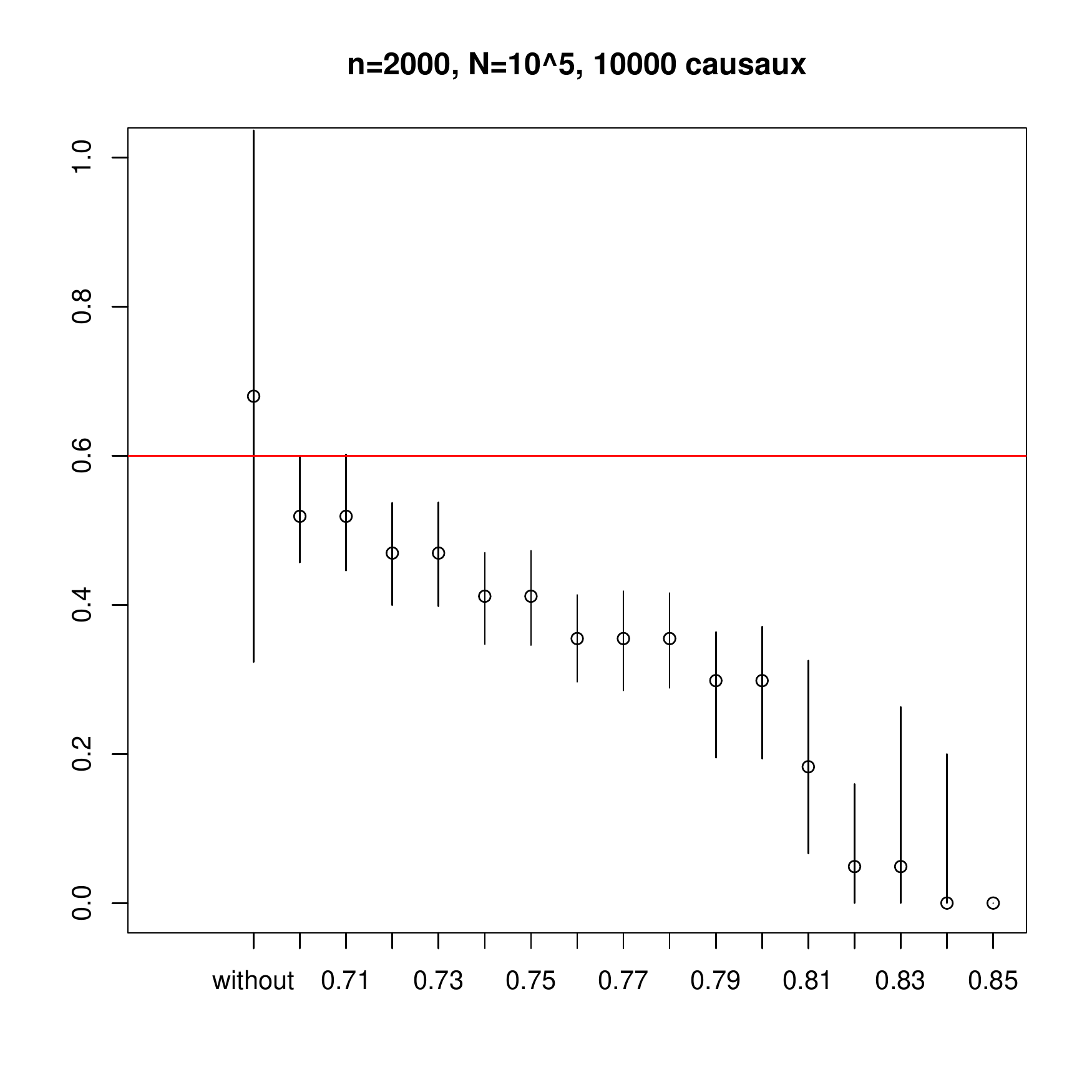}\\
\end{tabular}
\caption{\textcolor{black}{Estimation of the heritability with 95\% confidence intervals for $\eta^\star$ from 0.4 to 0.6 (from left to right), and from 100, 1000 and 10000 causal SNPs from top to bottom. Each graph shows the heritability estimations with 95\% confidence intervals computed with HiLMM (``without'') and for thresholds between 0.7 and 0.85.}}
  \label{fig:overlap}
\end{figure}

\textcolor{black}{On the one hand, we observed that applying HiLMM provides unbiased estimations of the heritability, no matter the number of causal SNPs. 
However, the main drawback of this estimator is its very large variance. On the other hand, if the number of causal SNPs is not too high, 
EstHer provides unbiased estimations of the heritability with standard errors substantially smaller than HiLMM. 
However, if the number of causal SNPs is high, EstHer underestimates
the heritability. \textcolor{black}{These observations are similar to
  those made by \cite{zhou:carbonetto:stephens:2013}, who built an
  hybrid estimator able to deal with both sparse and non sparse
  scenario, to which we will compare our approach in Section \ref{sec:comparison}.}
Therefore, we propose hereafter a rule to decide whether it is better to apply EstHer or HiLMM.
We can see from Figure \ref{fig:choix_seuil} that when there are 100 causal SNPs, there is a large range of threshold values which 
provide an accurate estimation of $\eta^\star$, but when there are 1000 or 10000 causal SNPs, see Figure \ref{fig:no_common_threshold}), 
the estimations are very different even for close thresholds. This observation gave us the idea of quantifying the stability of the estimations 
around the threshold that we determined as the optimal one. More precisely, for each threshold, we have an estimation of heritability with a 95\% confidence interval, 
and we count the number of thresholds for which the confidence intervals overlap. Figure \ref{fig:overlap} confirms the stability around the best threshold 
for different values of $\eta^\star$ and Table \ref{tab:table2} displays the number of ovelapping confidence intervals. 
We empirically determine the following criterion: if the mean number of thresholds is greater than 10 (over 16 tested thresholds), we apply EstHer, if not, we apply HiLMM. 
The results obtained by using this criterion are displayed in Figure \ref{compar_bslmm}.
}

 \begin{table}
   \caption{\textcolor{black}{Mean value of the number of overlapping confidence intervals for 16 thresholds from 0.7 to 0.85.}}
\begin{tabular}{|l| l |l |l|}
  \hline
 $ \eta^{\star} $ & 100 causal SNPs & 1000 causal SNPs & 10000 causal SNPs  \\
  \hline
  0.4 & 11.9 & 8.1 & 7.8 \\
  \hline   
  0.5 & 15.3 & 6.9 & 7  \\
  \hline
  0.6 & 16 & 9.2 & 7.1 \\
  \hline
\end{tabular}
\label{tab:table2}
\end{table}

\section{\textcolor{black}{Results after applying the decision criterion
    and comparison to other methods}}
\label{sec:comparison}

\textcolor{black}{
\subsection{Statistical performances}
In this section we show the results obtained after
  applying the criterion described in Section \ref{sec:criterion}. 
We compare these results to those obtained using HiLMM, but also with
the software GEMMA described in \cite{zhou:stephens:2012}. 
GEMMA can fit both a non sparse linear mixed model (GEMMA-LMM) and a sparse linear
mixed model if the BSLMM option is chosen denoted by BSLMM in the sequel. As explained in
\cite{zhou:carbonetto:stephens:2013}, BSLMM can deal with very sparse
and also with very polygenic scenarios. 
}

\begin{figure}[!ht]
  \centering
  \begin{tabular}{cc}
  
 \includegraphics[trim=0mm 15mm 0mm 20mm,clip,width=.45\textwidth] {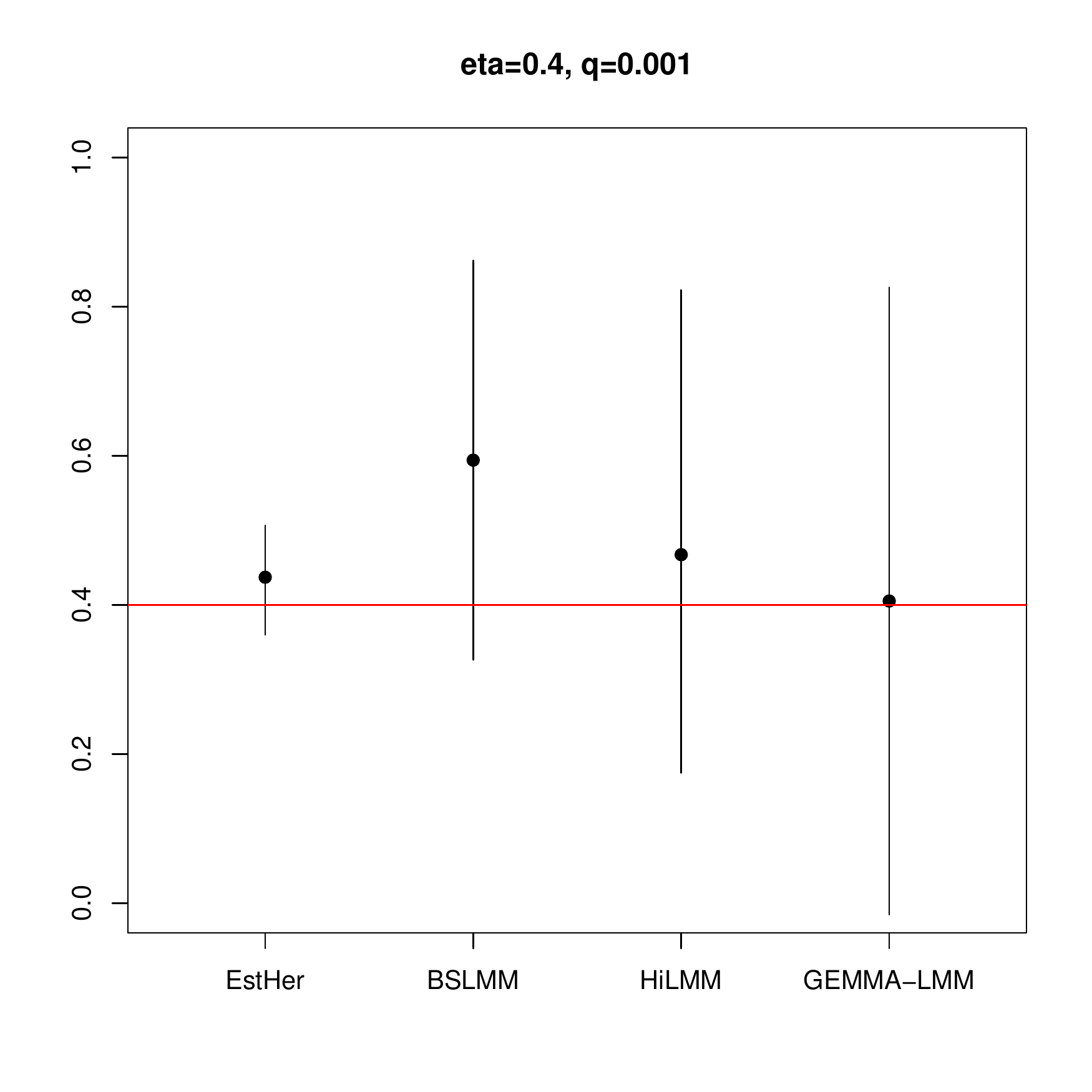} 
  & \includegraphics[trim=0mm 15mm 0mm 20mm,clip,width=.45\textwidth]{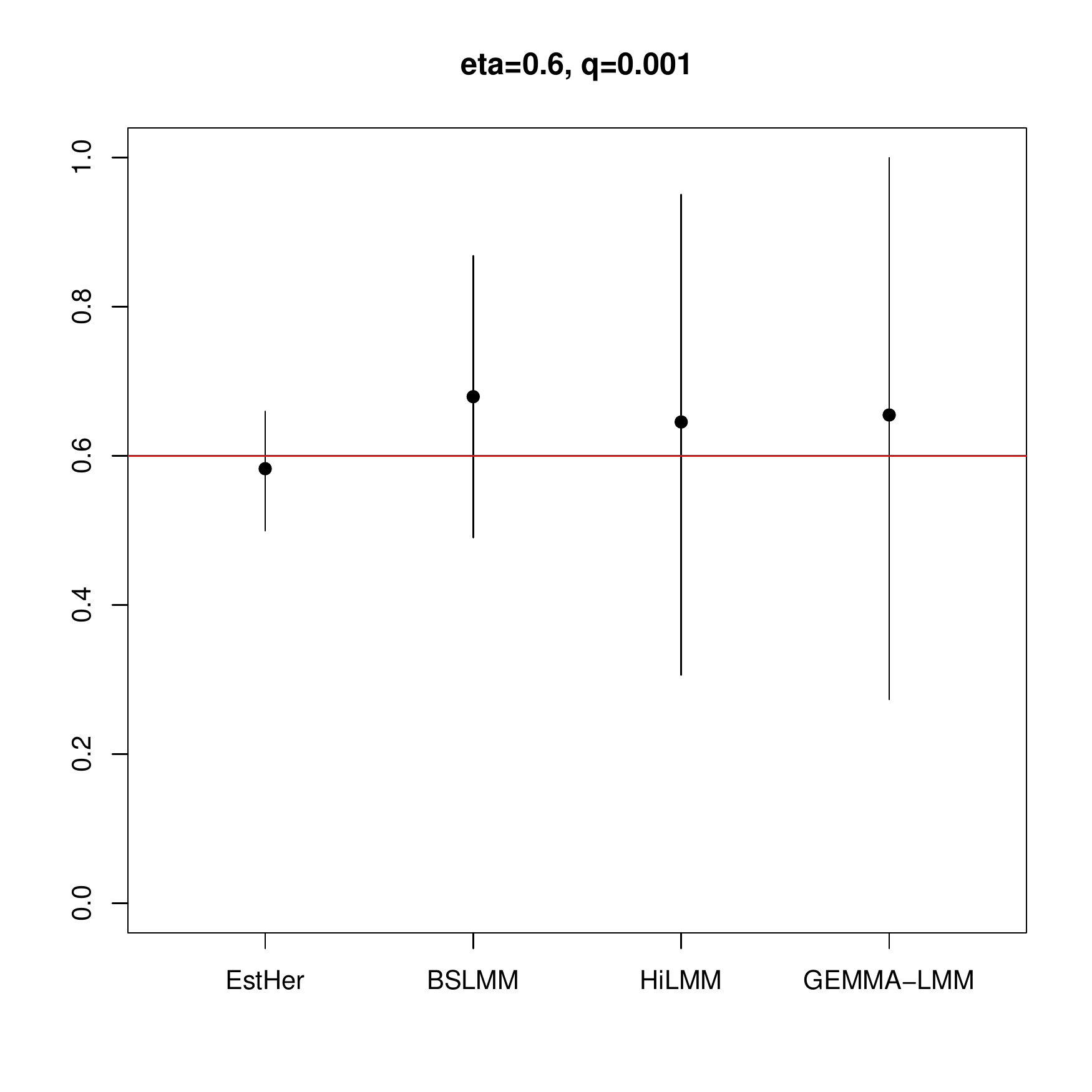}
\\
  \includegraphics[trim=0mm 15mm 0mm 20mm,clip,width=.45\textwidth]{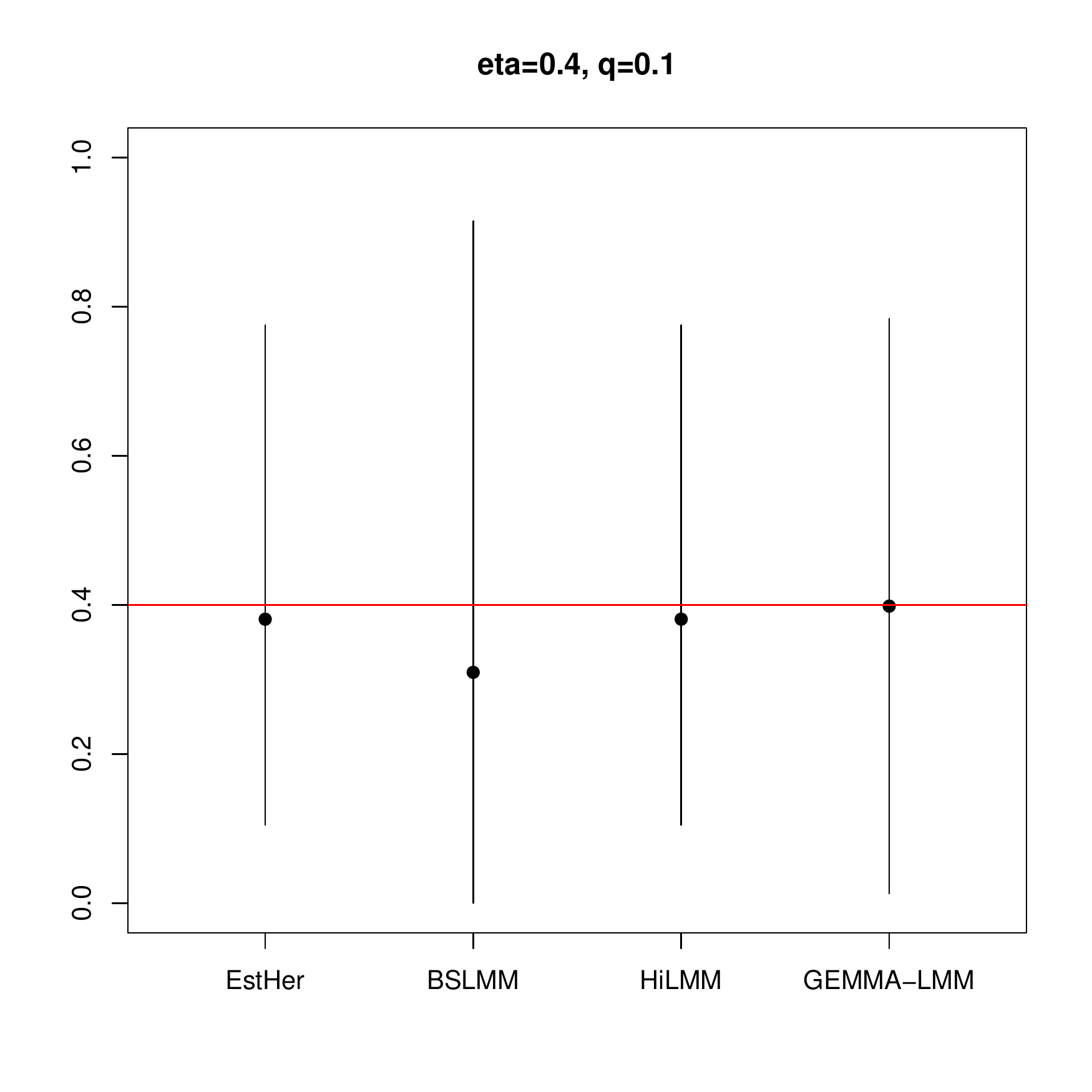}
  & \includegraphics[trim=0mm 15mm 0mm 20mm,clip,width=.45\textwidth]{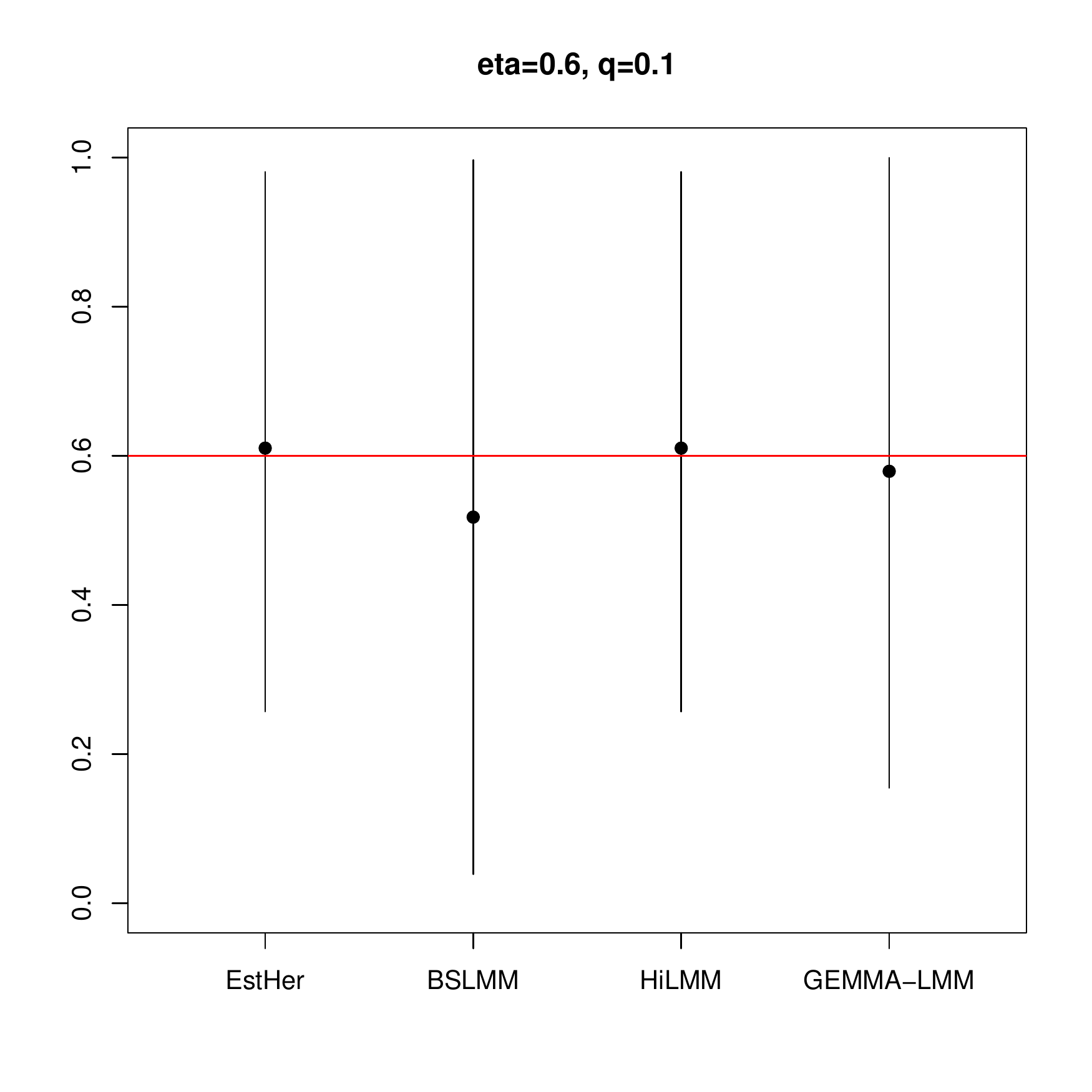}\\
  \end{tabular}
  \caption{Estimations of $\hat{\eta}$ with 95 \% confidence intervals
    obtained using EstHer, BSLMM, HiLMM and GEMMA-LMM with 100 causal SNPs
    (top) and 10,000 causal SNPs (bottom). The results are
    obtained with 10 replications.}
  \label{compar_bslmm}
\end{figure}

\textcolor{black}{
We can see from the bottom
part of Figure \ref{compar_bslmm} 
that, in very polygenic scenarios ($q=0.1$, namely 10,000 causal SNPs), all the methods provide similar
results: the four estimators are indeed empirically 
unbiased, but with a very large variance. 
}

\textcolor{black}{
In sparse scenarios ($q=10^{-3}$, namely 100 causal SNPs), we can see from the top part of Figure
\ref{compar_bslmm} that EstHer provides better results than HiLMM and
GEMMA-LMM which exhibit similar statistical performances. 
In sparse scenarios, the variance of the BSLMM estimator is larger
than the one provided by  EstHer and smaller than the one provided by
GEMMA-LMM and HiLMM. However, the performances of BSLMM could perhaps
be improved by changing the MCMC parameters. Here, for
computational time reasons, we used the default parameters that is 100,000 and 1,000,000
for the number of burn-in steps and the number of sampling, respectively.
}

\subsection{\textcolor{black}{Computational times}}
\begin{figure}[!ht]
\includegraphics[trim=0mm 15mm 0mm 20mm,clip,width=.5\textwidth]{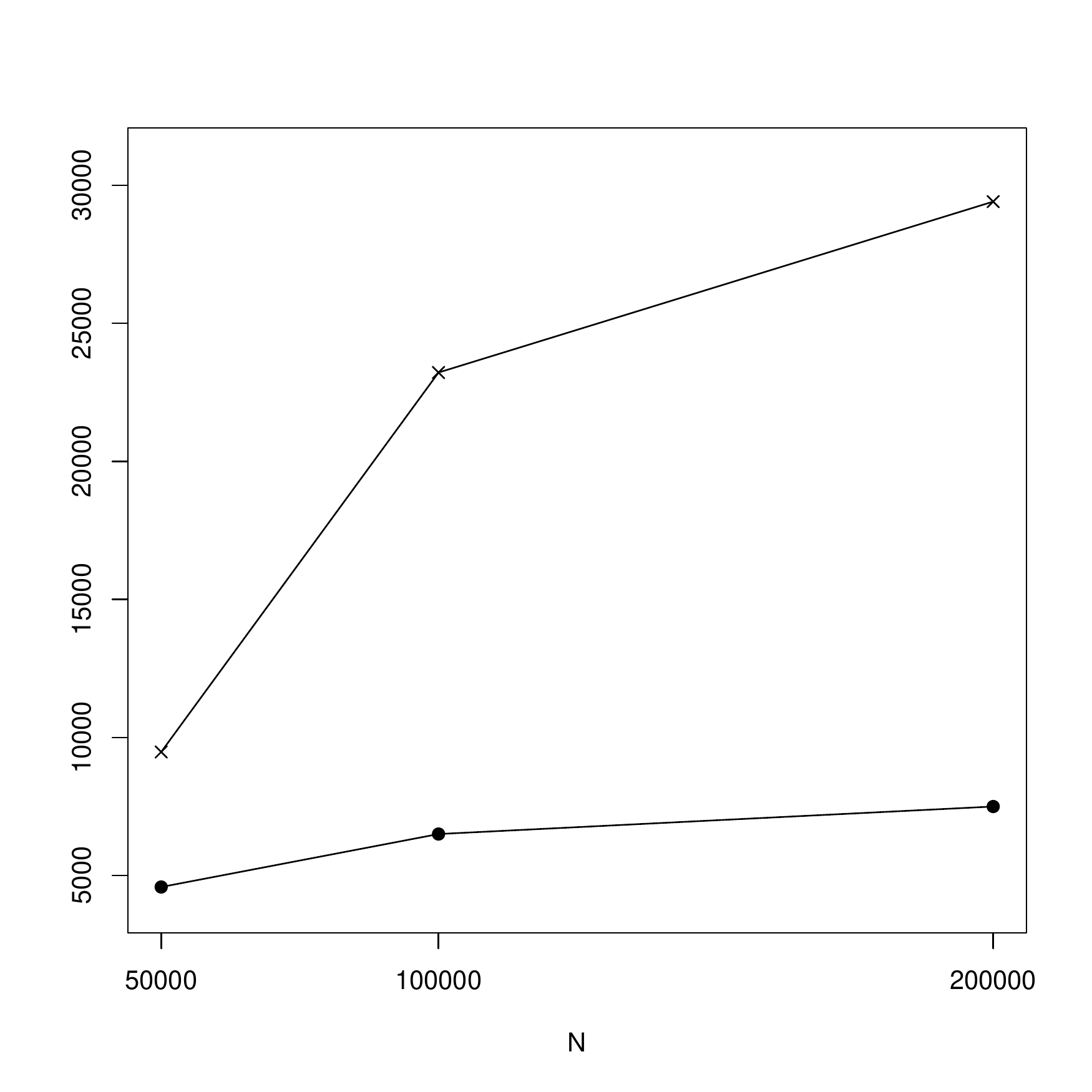}
 \caption{\textcolor{black}{Times (in seconds) to compute one heritability estimation with BSLMM (crosses)
and EstHer (dots) by using 16 thresholds for $n=2000$ and different values of $N$ from $50,000$ to $200,000$.}}
 \label{fig:compar_tps}
\end{figure}

\textcolor{black}{
The computational times in seconds for one estimation of the heritability with
BSLMM and the heritability estimation for 16 thresholds as well as the
associated confidence intervals with our method EstHer are displayed
in Figure \ref{fig:compar_tps}.
We chose this number of thresholds since we applied the criterion defined in Section \ref{sec:criterion}.
It should be noticed that  the computational times for EstHer could be
reduced by diminishing the number of thresholds. For BSLMM we used the
default parameters for the number of burn-in steps and the number of
sampling. We can see from this figure that the gap between EstHer and
BSLMM is all the more important that $N$ is large. Contrary to our
approach, BSLMM seems to be very sensitive in terms of computational
time to the value of $N$.
}


\section{Applications to genetic data}
\label{sec:real_data}

In this section, we applied our method to the neuroanatomic data coming from the Imagen project.
\textcolor{black}{In this data set, $n=2087$ individuals and $N=273926$ SNPs. For further details on this data set, we refer the reader to Section \ref{sec:data_set}.}

\subsection{Calibration of the threshold}

\begin{figure}[!ht]
  \centering
\begin{tabular}{cc}
  \includegraphics[trim=10mm 15mm 0mm 20mm,clip,width=.45\textwidth] {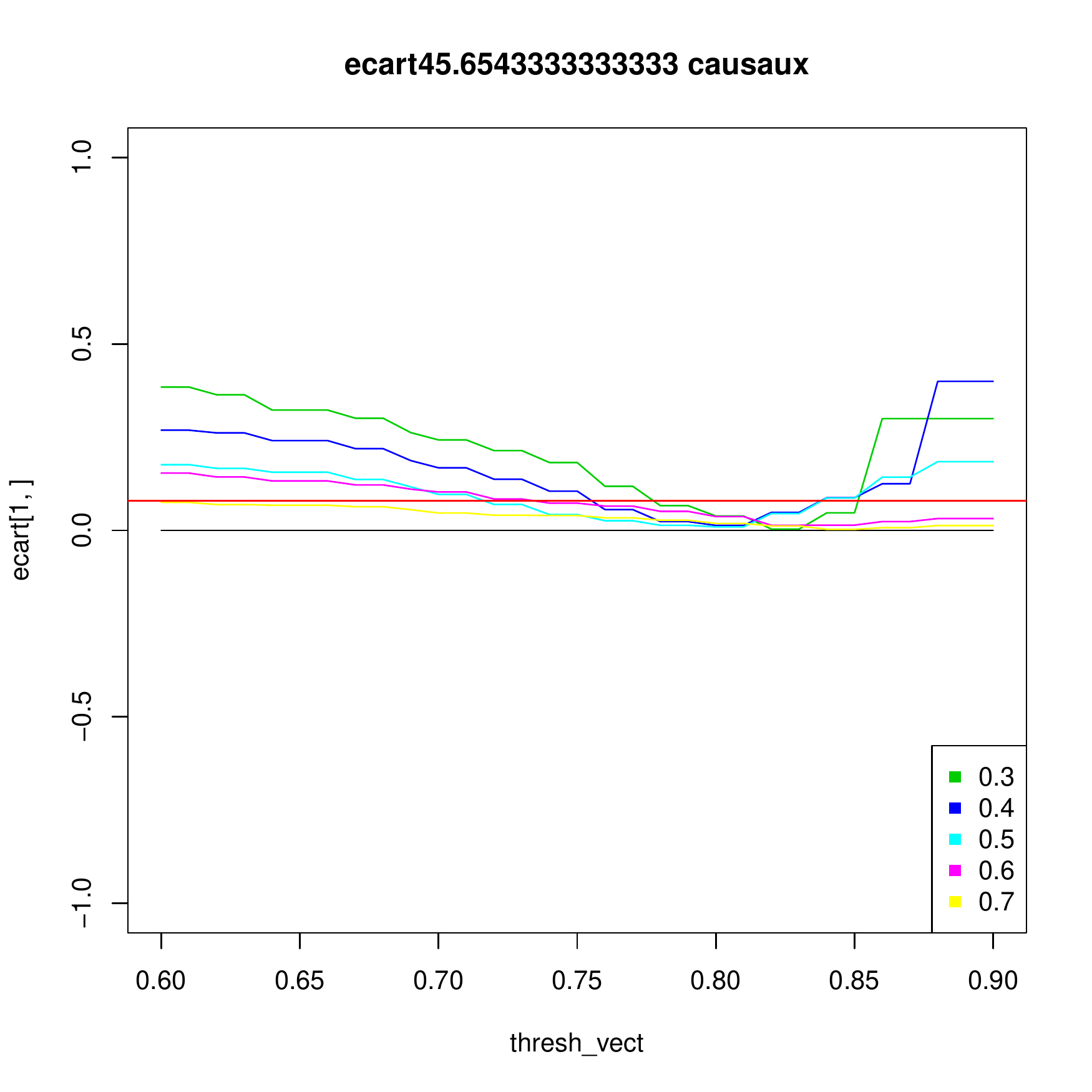}\label{Zvraie_50causaux}
&  \includegraphics[trim=10mm 15mm 0mm 20mm,clip,width=.45\textwidth]{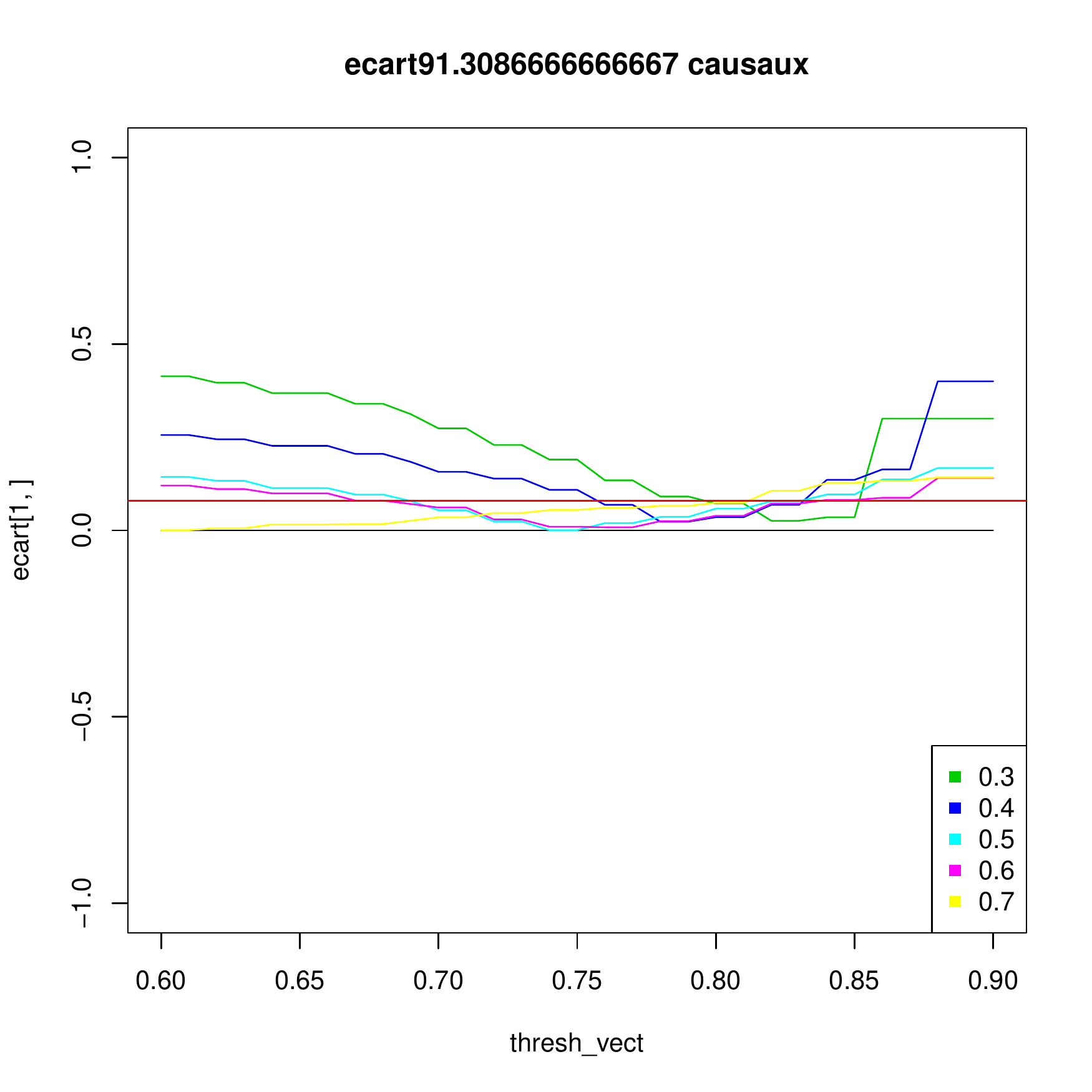}\label{Zvraie_100causaux} \\
(a) & (b)\\
\end{tabular}
  \caption{Absolute value of the difference between $\eta^\star$ and $\hat{\eta}$ for thresholds from 0.6 to 0.9, and for different values of $qN$: (a) 50 causal SNPs, 
(b)  100 causal SNPs. 
Each difference has been computed as the mean of 10 replications.}
  \label{fig:choix_seuil_ZVraie}
\end{figure}

We start by finding the threshold which is the most adapted to the Imagen data set. 
\textcolor{black}{We use the same technique as the one described in Section \ref{subsec:thresh_choice}: for several values of $\eta^\star$ 
and several thresholds, we display the absolute value of $\eta^\star-\hat{\eta}$, see Figure \ref{fig:choix_seuil_ZVraie}. 
The only difference with Section \ref{subsec:thresh_choice} is that we generated the observations by using the matrix 
$\W$ coming from the IMAGEN data set. According to Figure \ref{fig:choix_seuil_ZVraie}, we can find a
reliable range of thresholds for estimating the heritability for all $\eta^\star$ from 0.4 to 0.7 when the number of causal SNPs is smaller than 100.
This optimal threshold is equal to 0.79. We shall use this value in the sequel.}

\subsection{\textcolor{black}{Application of the decision criterion}}

\textcolor{black}{Since we determined in the previous section that the optimal threshold is 0.79, we apply EstHer for thresholds around this value, that is from 0.7 to 0.85. 
We then count the number of overlapping confidence intervals, as explained in Section \ref{sec:criterion}. The results are displayed in Table \ref{tab:table3}. 
We observe from this table that the sensitivity to the choice of the threshold varies substantially from one phenotype to another. 
Hence, we choose to apply our EstHer approach to the most stable phenotypes with respect to our criterion, namely pa, amy and acc. 
For the other phenotypes we recommand to apply HiLMM or another similar approach such as GCTA or GEMMA-LMM.}

 \begin{table}
   \caption{Mean value of the number of overlapping confidence intervals for 16 thresholds from 0.7 to 0.85.}
\begin{tabular}{|l| l |}
  \hline
Phenotype & Number of thresholds  \\
  \hline
  Bv & 7.19 \\
  \hline   
  Hip & 7.5  \\
  \hline
 Icv & 7.37 \\
 \hline 
 Acc & 9.94 \\
 \hline
 Amy & 9.88 \\
  \hline
  Th & 7.5 \\
  \hline
  Ca & 7.13 \\
  \hline
  Pu & 7.13 \\
  \hline
  Pa & 10.75 \\
  \hline
\end{tabular}
\label{tab:table3}
\end{table}

\subsection{Results}

\begin{figure}[!ht]
  \centering
\begin{tabular}{cc}
\hspace{-5mm}\includegraphics[scale=0.35]{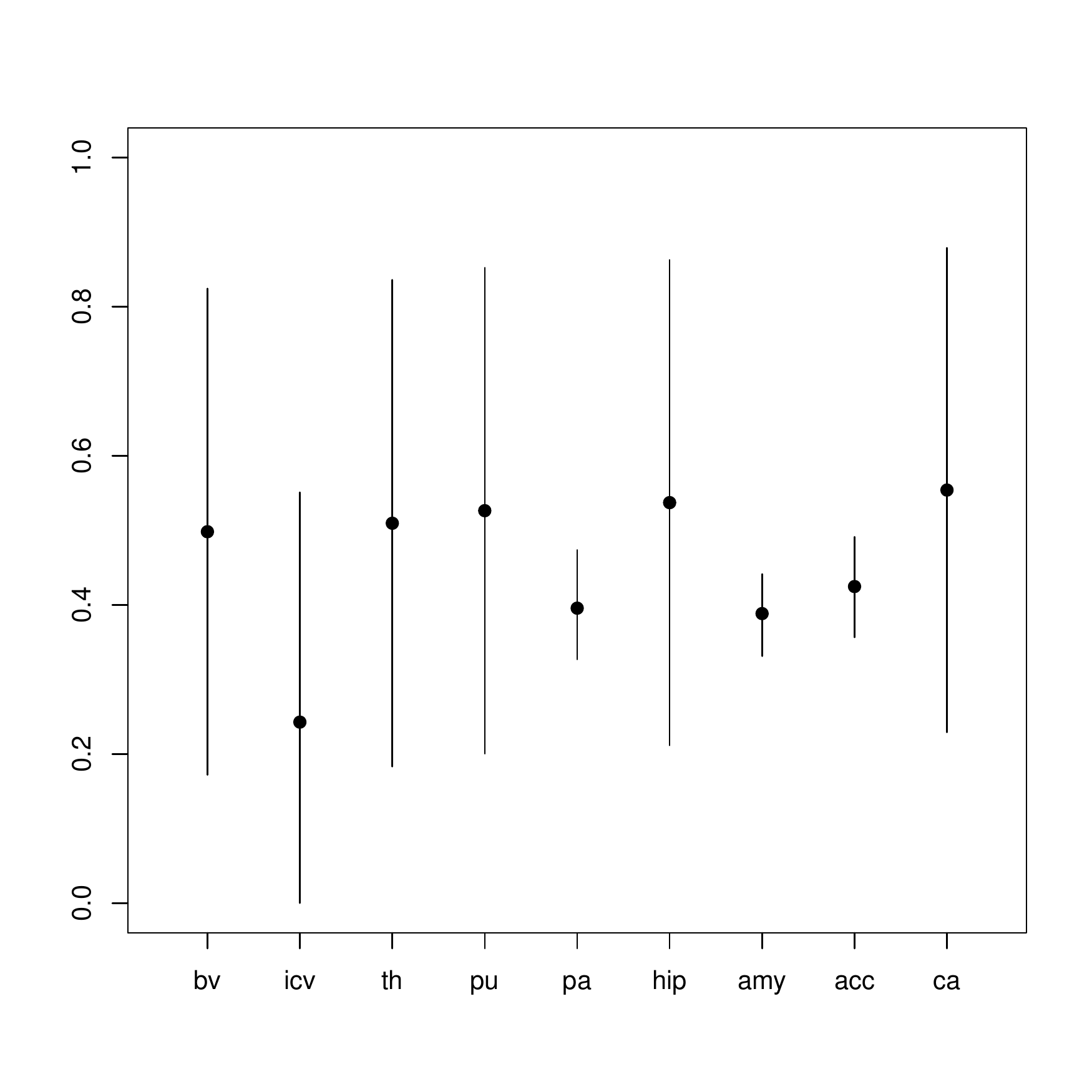}
 & \includegraphics[scale=0.35]{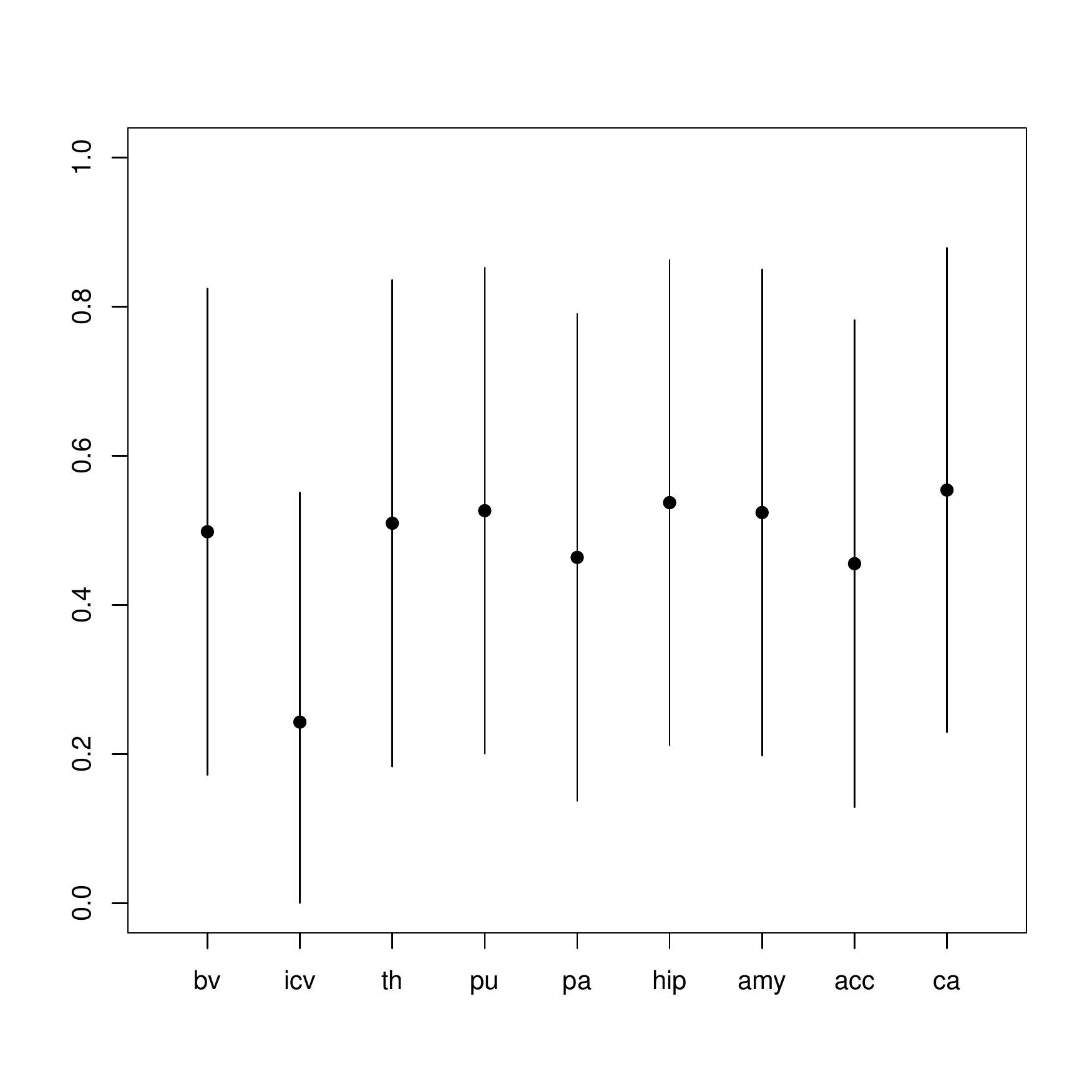} \\
 (a) & (b) \\
 \end{tabular}
  \caption{(a) Heritability estimations of bv, icv, th, pu, pa, hip,
    amy, acc, and ca with 95\% confidence intervals obtained using
    EstHer or HiLMM according to the outcome of our decision
    criterion. 
(b) Heritability estimations of bv, icv, th, pu, pa, hip, amy, acc and ca with 95\% confidence intervals obtained using HiLMM.}
  \label{fig:imagen}
\end{figure}

\textcolor{black}{Figure \ref{fig:imagen} (a) shows the heritability
  estimation with 95 \% confidence intervals for all phenotypes, using
  either EstHer or HiLMM according to 
the outcome of our decision criterion. Figure \ref{fig:imagen} (b)
shows the results obtained by using HiLMM, namely without any variable
selection step.} 
We compare our results with the ones obtained by \cite{toro:2014} who
estimated the heritability of the same phenotypes by using the software GCTA.
On the one hand, we can see from Figure \ref{fig:imagen} that in the cases where EstHer
is used the confidence intervals given by our methodology are
substantially smaller and included in those provided by either HiLMM
or \cite{toro:2014}. On the
other hand, when HiLMM is used our results are on a par with those
obtained by \cite{toro:2014}.
Moreover, our approach provides a list of SNPs which may contribute to the variations of a given 
phenotype and which could be further analyzed from a biological point of view in order to identify new biological pathways.



\section{Conclusion}

We propose in this paper a practical method to estimate the heritability in sparse linear mixed models 
using variable selection tools, as well as confidence interval obtained thanks to a  non parametric bootstrap approach.
Our approach is implemented in the R package EstHer which is available from the Comprehensive R Archive Network (CRAN) and from the web 
page of the first author. 
In the course of this study, we showed that our approach has two main features which makes it very attractive.
Firstly, it is very efficient from a statistical point of view since it provides confidence intervals considerably smaller than those obtained with 
methods without variable selection. Secondly, its very low computational burden makes its use
feasible on very large data sets coming from quantitative genetics.

\textcolor{black}{Moreover, we observed that the statistical performance of the EstHer approach
are all the more impressive that the level of sparsity is high that is when $q$ is small. For this reason,
we also proposed an empirical criterion which allows the user
to decide whether it is better to apply an approach that takes into account the sparsity and starts with a variable selection stage, namely EstHer,
or an approach which ignores the potential sparsity in the observations, namely HiLMM.}

\section*{Acknowledgments}

The authors would like to thank Nicolai Meinshausen and Nicolas Verzelen for fruitful discussions and the IMAGEN consortium for providing the data.

\bibliographystyle{abbrv}
\bibliography{bibli_pap2}

%
%
%
%
%
%
%

\end{document}